\documentclass[twoside,12pt,a4wide]{article}
\usepackage{amssymb}

\usepackage{makeidx}
\usepackage{a4wide}
\usepackage{amsthm}
\usepackage{amsmath}
\usepackage[all]{xy}
\usepackage{enumerate}
\usepackage{fancyhdr}

\newtheorem{theorem}{Theorem}[section]
\newtheorem{proposition}[theorem]{Proposition}
\newtheorem{corollary}[theorem]{Corollary}
\newtheorem{lemma}[theorem]{Lemma}
\theoremstyle{definition}

\newtheorem{c-example}[theorem]{Counter Example}

\newtheorem*{Beweis}{Proof}
\newtheorem{definition}[theorem]{Definition}
\newtheorem{punto}[theorem]{}
\theoremstyle{remark}
\newtheorem{remark}[theorem]{Remark}

\CompileMatrices
\input{tcilatex}

\begin{document}

\title{A Note On Coinduction Functors between Categories of Comodules for Corings%
\thanks{%
MSC (2000): 16W30, 16D90, 18A40 \newline
Keywords: Corings, Comodules, Locally Projective Modules, Coinduction
Functors, Cotensor Product}}
\author{\textbf{Jawad Y. Abuhlail}\thanks{%
Supported by KFUPM} \\
Department of Mathematical Sciences\\
King Fahd University of Petroleum $\&\;$Minerals\\
31261 Dhahran - Saudi Arabia \\
abuhlail@kfupm.edu.sa}
\date{}
\maketitle

\begin{abstract}
In this note we consider different versions of coinduction functors between
categories of comodules for corings induced by a morphism of corings. In
particular we introduce a new version of the coinduction functor in the case
of \emph{locally projective} corings as a composition of suitable ``Trace''
and ``Hom'' functors and show how to derive it from a more \emph{general }%
coinduction functor between categories of type $\sigma \lbrack M].$ In
special cases (e.g. the corings morphism is part of a morphism of \emph{%
measuring }$\alpha $\emph{-pairings }or the corings have the same base
ring), a version of our functor is shown to be isomorphic to the usual
coinduction functor obtained by means of the cotensor product. Our results
in this note generalize previous results of the author on coinduction
functors between categories of comodules for coalgebras over commutative
base rings.
\end{abstract}

\section{Introduction}

With $R$ we denote a commutative ring with $1_{R}\neq 0_{R},$ with $A$ an
arbitrary $R$-algebra and with $\mathbb{M}_{A}$ (respectively $_{A}\mathbb{M}
$) the category of right (respectively left) $A$-modules. The unadorned $%
-\otimes -,$ $\mathrm{Hom}(-,-)$ mean $-\otimes _{R}-,$ $\mathrm{Hom}%
_{R}(-,-)$ respectively. For a right (respectively a left) $A$-module $M$ we
denote by 
\begin{equation*}
\vartheta _{M}^{r}:M\otimes _{A}A\rightarrow M\text{ (respectively }%
\vartheta _{M}^{l}:A\otimes _{A}M\rightarrow M\text{)}
\end{equation*}
the canonical $A$-isomorphism.

By{\normalsize \ }an $A$\emph{-ring}\ we mean an $A$-bimodule $\mathcal{T}$
with a $A$-bilinear maps $\mu _{T}:\mathcal{T}\otimes _{A}\mathcal{T}%
\rightarrow \mathcal{T}$ and $\eta _{T}:A\rightarrow \mathcal{T},$ such that 
$\mu _{\mathcal{T}}\circ (\mu _{\mathcal{T}}\otimes _{A}id_{\mathcal{T}%
})=\mu _{\mathcal{T}}\circ (id_{\mathcal{T}}\otimes _{A}\mu _{\mathcal{T}}),$
$\mu _{\mathcal{T}}\circ (id_{\mathcal{T}}\otimes _{A}\eta _{\mathcal{T}%
})=\vartheta _{\mathcal{T}}^{r}$ and $\mu _{\mathcal{T}}\circ (\eta _{%
\mathcal{T}}\otimes _{A}id_{\mathcal{T}})=\vartheta _{\mathcal{T}}^{l}.$ If $%
\mathcal{T}$ and $\mathcal{S}$ are $A$-rings with unities $\eta _{\mathcal{T}%
},$ $\eta _{\mathcal{S}}$ respectively, then an $A$-bilinear map $f:\mathcal{%
T}\rightarrow \mathcal{S}$ is called a \emph{morphism of }$A$\emph{-rings},
if $f\circ \mu _{\mathcal{T}}=\mu _{\mathcal{S}}\circ (f\otimes _{A}f)$ and $%
f\circ \eta _{\mathcal{T}}=\eta _{\mathcal{S}}.$ For an $A$-ring $\mathcal{T}
$ and two left (respectively right) $\mathcal{T}$-modules $M$ and $N$ we
denote with $\mathrm{Hom}_{\mathcal{T}-}(M,N)$ (respectively $\mathrm{Hom}_{-%
\mathcal{T}}(M,N)$) the set of all $\mathcal{T}$-linear maps from $M$ to $N.$
If $M$ and $N$ are $\mathcal{T}$-bimodules, then $\mathrm{Hom}_{\mathcal{T}-%
\mathcal{T}}(M,N)$ denotes the set of all $\mathcal{T}$-bilinear maps from $M
$ to $N.$

Let $M$ and $N$ be right $A$-modules and $f:M\rightarrow N$ be $A$-linear.
For a left $A$-module $L,$ the morphism $f$ will be called $L$\emph{-pure},
if the following sequence is exact 
\begin{equation*}
0\longrightarrow \mathrm{Ker}(f)\otimes _{A}L\overset{\iota \otimes
_{A}id_{L}}{\longrightarrow }M\otimes _{A}L\overset{f\otimes _{A}id_{L}}{%
\longrightarrow }N\otimes L.
\end{equation*}
If $f$ is $L$-pure for all left $A$-modules $L$ then $f$ is called a \emph{%
pure morphism }(e.g. \cite[40.13]{BW03}). If $M\subseteq N$ is a right $A$%
-submodule then it is called a \emph{pure submodule, }provided that the
embedding $M\overset{\iota }{\hookrightarrow }N$ is a pure morphism
(equivalently, if $\iota \otimes _{A}id_{L}:M\otimes _{A}L\rightarrow
N\otimes _{A}L$ is injective for every left $A$-module $L$). Pure morphisms
and pure submodules in the category of left $A$-modules are defined
analogously. A morphism of $A$-bimodules is said to be pure, if it is pure
in $\mathbb{M}_{A}$ as well as in $_{A}\mathbb{M}.$ For an $A$-bimodule $N,$
we call an $A$-subbimodule $M\subset N$ \emph{pure}, if $M_{A}\subset N_{A}$
and $_{A}M\subset $ $_{A}N$ are pure.

A left (respectively right) $A$-module $W$ is called \emph{locally projective%
} (in the sense of B. Zimmermann-Huisgen \cite{Z-H76}), if for every diagram
of left (respectively right) $A$-modules 
\begin{equation*}
\xymatrix{0 \ar[r] & F \ar@{.>}[dr]_{g' \circ \iota} \ar[r]^{\iota} & W
\ar[dr]^{g} \ar@{.>}[d]^{g'} & & \\ & & L \ar[r]_{\pi} & N \ar[r] & 0}
\end{equation*}
with exact rows and $F$ f.g.: for every $A$-linear map $g:W\rightarrow N,$
there exists an $A$-linear map $g^{\prime }:W\rightarrow L,$ such that $%
g\circ \iota =\pi \circ g^{\prime }\circ \iota .$ Note that every projective
left (respectively right) $A$-module is locally projective.

Let $A,$ $B$ be $R$-algebras with an $R$-algebra morphism $\beta
:A\rightarrow B$ and $M$ (respectively $N$) be a right (respectively a left) 
$B$-module. We consider $M$ (respectively $N$) as a right (respectively as a
left ) $A$-module through 
\begin{equation*}
m\leftharpoonup a:=m\beta (a)\text{ (respectively }a\rightharpoonup n:=\beta
(a)n\text{)}
\end{equation*}
and denote with 
\begin{equation*}
\chi _{M,N}:M\otimes _{A}N\rightarrow M\otimes _{B}N
\end{equation*}
the canonical $R$-linear morphism.

\bigskip

(Co)induction functors between categories of entwined modules induced by a
morphism of entwined structures were studied in the case of base fields by
T. Brzezi\'{n}ski et al. (e.g. \cite{Brz99}, \cite{BCMZ01}). In the general
case of a morphism of corings $(\theta :\beta ):(\mathcal{C}:A)\rightarrow (%
\mathcal{D}:B),$ J. G\'{o}mez-Torrecillas presented in \cite{G-T02} (under
some purity conditions and using the cotensor product) a coinduction functor 
$\mathcal{G}:\mathbb{M}^{\mathcal{D}}\rightarrow \mathbb{M}^{\mathcal{C}}$
and proved it is right adjoint to the canonical induction functor $-\otimes
_{A}B:\mathbb{M}^{\mathcal{C}}\rightarrow \mathbb{M}^{\mathcal{D}}$ ( 
\cite[Proposition 5.3]{G-T02}).

\bigskip

In his dissertation \cite{Abu01}, the author presented and studied
coinduction functors between categories of comodules induced by morphisms
between \emph{measuring }$\alpha $\emph{-pairings }for coalgebras (over a
common commutative base rings). In this note the author generalizes his
previous results to the case of coinduction functors between categories of
comodules for corings over (possibly different) arbitrary base rings.
Although in our setting some technical assumptions are assumed and the
corings under consideration are restricted to the locally projective ones,
the main advantages of the new version of the coinduction functor is that it
is presented as a composition of well-studied ``\textrm{Trace}'' and ``%
\textrm{Hom}'' functors and that it avoids the use of the cotensor product.
Another advantage of the new version of the coinduction functor we introduce
is that it is derived from a more general coinduction functor between
categories of type $\sigma \lbrack M],$ the theory of which is well
developed (e.g. \cite{Wis88}, \cite{Wis96}). With these two main advantages
in mind, this note aims to enable a more intensive study of the coinduction
functors between categories of comodules and their properties which will be
carried out elsewhere.

After this introductory section, we present in the second section some
definitions and lemmas, that will be used later. In the third section we
consider coinduction functors between categories of type $\sigma \lbrack M]$
as a general model from which we derive later our version of coinduction
functors between categories of comodules for locally projective corings. In
the fourth section we consider the case of coinduction functors induced by
morphisms of measuring $\alpha $-pairings, which turn out to be special
forms of the coinduction functors introduced in the third section. We also
handle the special case, where the corings are defined over a common base
ring. In particular we prove that our coinduction functor in these cases is
isomorphic to the coinduction functor presented by J. G\'{o}mez-Torrecillas
in \cite{G-T02} by means of the cotensor product. In the fifth and last
section we handle the general case of coinduction functors induced by
morphisms of corings (over possibly different) arbitrary base rings.

\section{Preliminaries}

In this section we present the needed definitions and lemmas.

\begin{punto}
\textbf{Subgenerators}.\label{subg} Let $A$ be an $R$-algebra, $\mathcal{T}$
be an $A$-ring and $K$ be a right $\mathcal{T}$-module. We say a right $%
\mathcal{T}$-module $N$ is $K$\emph{-subgenerated} if $N$ is isomorphic to a
submodule of a $K$-generated right $\mathcal{T}$-module, equivalently if $N$
is the kernel of a morphism between $K$-generated right $\mathcal{T}$%
-modules. With $\sigma \lbrack K_{\mathcal{T}}]\subseteq \mathbb{M}_{%
\mathcal{T}}$ we denote the \emph{full }subcategory of $\mathbb{M}_{\mathcal{%
T}}$ whose objects are $K$-subgenerated. In fact $\sigma \lbrack K_{\mathcal{%
T}}]\subseteq \mathbb{M}_{\mathcal{T}}$ is the \emph{smallest }Grothendieck
full subcategory that contains $K.$ Moreover we have the \emph{trace functor}
\begin{equation*}
\mathrm{Sp}(\sigma \lbrack K_{\mathcal{T}}],-):\mathbb{M}_{\mathcal{T}%
}\rightarrow \sigma \lbrack K_{\mathcal{T}}],\text{ }M\mapsto \sum \{f(N):%
\text{ }f\in \mathrm{Hom}_{-\mathcal{T}}(N,M),\text{ }N\in \sigma \lbrack K_{%
\mathcal{T}}]\},
\end{equation*}
which is right-adjoint to the inclusion functor $\iota :\sigma \lbrack K_{%
\mathcal{T}}]\hookrightarrow \mathbb{M}_{\mathcal{T}}$ (\cite[45.11]{Wis88}%
). In fact $\mathrm{Sp}(\sigma \lbrack K_{\mathcal{T}}],M)$ is the largest
right $\mathcal{T}$-submodule of $M_{\mathcal{T}}$ that belongs to $\sigma
\lbrack K_{\mathcal{T}}]$ and $m\in \mathrm{Sp}(\sigma \lbrack K_{\mathcal{T}%
}],M)$ if and only if there exists a finite subset $W\subset K$ with $%
\mathrm{Ann}_{-A}(W)\subseteq \mathrm{Ann}_{-A}(m).$ For a left $\mathcal{T}$%
-module $L,$ the full category $\sigma \lbrack _{\mathcal{T}}L]\subseteq $ $%
_{\mathcal{T}}\mathbb{M}$ is defined analogously. The reader is referred to 
\cite{Wis88} and \cite{Wis96} for the well developed theory of categories of
this type.
\end{punto}

\begin{definition}
Let $A$ be an $R$-algebra. A \emph{coassociative} $A$\emph{-coring }$(%
\mathcal{C},\Delta _{\mathcal{C}},\varepsilon _{\mathcal{C}})$ is an $A$%
-bimodule $\mathcal{C}$ with $A$-bilinear maps 
\begin{equation*}
\Delta _{\mathcal{C}}:\mathcal{C}\rightarrow \mathcal{C}\otimes _{A}\mathcal{%
C},\text{ }c\mapsto \sum c_{1}\otimes _{A}c_{2}\text{ and }\varepsilon _{%
\mathcal{C}}:\mathcal{C}\rightarrow A,
\end{equation*}
such that the following diagrams are commutative 
\begin{equation*}
\begin{tabular}{lll}
$\xymatrix{ {\mathcal {C}} \ar^(.45){\Delta_{\mathcal C}}[rr]
\ar_(.45){\Delta_{\mathcal C}}[d] & & {\mathcal {C}} \otimes_{A} {\mathcal
{C}} \ar^(.45){id_{\mathcal C} \otimes_A \Delta_{\mathcal C}}[d]\\ {\mathcal
{C}} \otimes_{A} {\mathcal {C}} \ar_(.45){\Delta_{\mathcal C} \otimes_A
id_{\mathcal C}}[rr] & & {\mathcal {C}} \otimes_{A} {\mathcal {C}}
\otimes_{A} {\mathcal {C}} }$ &  & $\xymatrix{ & & {\mathcal {C}}
\ar^(.45){\Delta _{\mathcal {C}}}[d] & & \\ {A} \otimes_{A} {\mathcal {C}}
\ar[urr]^(.45){\vartheta _{\mathcal {C}} ^l} & & {\mathcal {C}} \otimes_{A}
{\mathcal {C}} \ar^(.45){\varepsilon _{\mathcal {C}} \otimes_A id_{\mathcal
C}}[ll] \ar_(.45){id_{\mathcal C} \otimes_A \varepsilon _{\mathcal {C}}}[rr]
& & {\mathcal {C}} \otimes_{A} {A} \ar[ull]_(.45){\vartheta _{\mathcal {C}}
^r} }$%
\end{tabular}
\end{equation*}
The $A$-bilinear maps $\Delta _{\mathcal{C}},$ $\varepsilon _{\mathcal{C}}$
are called the \emph{comultiplication }and\emph{\ }the\emph{\ counit} of $%
\mathcal{C},$ respectively. All corings under consideration are assumed to
be coassociative and have counits. The objects in the \emph{category of
corings} $\mathbf{Crg}$ of coassociative corings are understood to be pairs $%
(\mathcal{C}:A),$ where $A$ is an $R$-algebra and $\mathcal{C}$ is an $A$%
\emph{-coring}. A \emph{morphism between corings} $(\mathcal{C}:A)$ and $(%
\mathcal{D}:B)$ consists of a pair of mappings $(\theta :\beta ):(\mathcal{C}%
:A)\rightarrow (\mathcal{D}:B),$ where $\beta :A\rightarrow B$ is an $R$%
-algebra morphism (so one can consider $\mathcal{D}$ as an $A$-bimodule) and 
$\theta :\mathcal{C}\rightarrow \mathcal{D}$ is an $(A,A)$-bilinear map with 
\begin{equation}
\chi _{\mathcal{D},\mathcal{D}}\circ (\theta \otimes _{A}\theta )\circ
\Delta _{\mathcal{C}}=\Delta _{\mathcal{D}}\circ \theta \text{ and }%
\varepsilon _{\mathcal{D}}\circ \theta =\beta \circ \varepsilon _{\mathcal{C}%
}.  \label{coring-mor}
\end{equation}
A morphism of corings $(\theta :\beta ):(\mathcal{C}:A)\rightarrow (\mathcal{%
D}:B)$ is said to be \emph{pure} if for every right $\mathcal{D}$-comodule $%
N,$ the following right $A$-module morphism is $\mathcal{C}$-pure (e.g. $_{A}%
\mathcal{C}$ is flat): 
\begin{equation*}
\varrho _{N}\otimes _{A}id_{\mathcal{C}}-(\chi _{N,\mathcal{D}}\otimes
_{A}id_{\mathcal{C}})\circ (id_{N}\otimes _{A}\theta \otimes id_{\mathcal{C}%
})\circ (id_{N}\otimes _{A}\Delta _{\mathcal{C}}):N\otimes _{A}\mathcal{C}%
\rightarrow N\otimes _{B}\mathcal{D}\otimes _{A}\mathcal{C}.
\end{equation*}
If $A=B$ and $\beta $ is not mentioned explicitly, then we mean $\beta
=id_{A}.$
\end{definition}

\begin{definition}
Let $(\mathcal{C},\Delta _{\mathcal{C}},\varepsilon _{\mathcal{C}})$ and $(%
\mathcal{D},\Delta _{\mathcal{D}},\varepsilon _{\mathcal{D}})$ be $A$\emph{%
-corings} with $_{A}\mathcal{D}_{A}\subseteq $ $_{A}\mathcal{C}_{A}.$ We
call $\mathcal{D}$ an $A$-subcoring of $\mathcal{C},$ if the embedding $%
\iota _{\mathcal{D}}:\mathcal{D}\rightarrow \mathcal{C}$ is a morphism of $A$%
-corings. If $_{A}\mathcal{D}\subseteq $ $_{A}\mathcal{C}$ and $\mathcal{D}%
_{A}\subseteq \mathcal{C}_{A}$ are pure, then $\mathcal{D}$ is a subcoring
of $\mathcal{C}$ if and only if $\Delta _{\mathcal{C}}(\mathcal{D})\subseteq 
\mathcal{D}\otimes _{A}\mathcal{D}.$
\end{definition}

\begin{punto}
\textbf{The dual rings of a coring. }(\cite{Guz85}) Let $(\mathcal{C},\Delta
_{\mathcal{C}},\varepsilon _{\mathcal{C}})$ be an $A$-coring. Then $^{\ast }%
\mathcal{C}:=\mathrm{Hom}_{A-}(\mathcal{C},A)$ is an associative $A$-ring
with multiplication given by the \emph{left convolution product} 
\begin{equation*}
(f\star _{l}g)(c):=\sum g(c_{1}f(c_{2}))\text{ for all }f,g\in \text{ }%
^{\ast }\mathcal{C},\text{ }c\in \mathcal{C},
\end{equation*}
unit $\varepsilon _{\mathcal{C}},$ and $\mathcal{C}^{\ast }:=\mathrm{Hom}%
_{-A}(\mathcal{C},A)$ is an associative $A$-ring with multiplication given
by the \emph{right convolution product} 
\begin{equation*}
(f\star _{r}g)(c):=\sum f(g(c_{1})c_{2})\text{ for all }f,g\in \text{ }%
\mathcal{C}^{\ast },\text{ }c\in \mathcal{C},
\end{equation*}
unit $\varepsilon _{\mathcal{C}}.$ Moreover $^{\ast }\mathcal{C}^{\ast }:=%
\mathrm{Hom}_{A-A}(\mathcal{C},A)$ is an associative $A$-ring with
multiplication given by the \emph{convolution product} 
\begin{equation*}
(f\star g)(c):=\sum g(c_{1})f(c_{2})\text{ for all }f,g\in \text{ }^{\ast }%
\mathcal{C}^{\ast },\text{ }c\in \mathcal{C}.
\end{equation*}
\end{punto}

\begin{punto}
Let $(\mathcal{C},\Delta _{\mathcal{C}},\varepsilon _{\mathcal{C}})$ be an $A
$-coring. A right $\mathcal{C}$\emph{-comodule} is a right $A$-module $M$
associated with a right $A$-linear map ($\mathcal{C}$\emph{-coaction}) 
\begin{equation*}
\varrho _{M}:M\rightarrow M\otimes _{A}\mathcal{C},\text{ }m\mapsto \sum
m_{<0>}\otimes _{A}m_{<1>},
\end{equation*}
such that the following diagram is commutative 
\begin{equation*}
\xymatrix{M \ar^(.4){\varrho _M}[rr] \ar_(.45){\varrho _M}[d] & & M
\otimes_{A} {\mathcal{C}} \ar^(.45){id _M \otimes_A \Delta_{\mathcal C}}[d]
\\ M \otimes_{A} {\mathcal{C}} \ar_(.4){\varrho _M \otimes_A
id_{\mathcal{C}}}[rr] & & M \otimes_{A} {\mathcal {C}} \otimes_{A}
{\mathcal{C}}}
\end{equation*}
We call $M$ \emph{counital}, if 
\begin{equation*}
\vartheta _{M}^{r}\circ (id_{M}\otimes \varepsilon _{\mathcal{C}})\circ
\varrho _{M}=id_{M}.
\end{equation*}
For right $\mathcal{C}$-comodules $(M,\varrho _{M}),$ $(N,\varrho _{N})\in 
\mathbb{M}^{\mathcal{C}}${\normalsize \ }we call a right $A$-linear map $%
f:M\rightarrow N$ a $\mathcal{C}$\emph{-comodule morphism}{\normalsize \ }%
(or \emph{right }$\mathcal{C}$\emph{-colinear}), if the following diagram is
commutative 
\begin{equation*}
\xymatrix{M \ar[rr]^{f} \ar[d]_{\varrho_M} & & N \ar[d]^{\varrho_N}\\ M
\otimes_{A} {\mathcal{C}} \ar[rr]_{f \otimes_A id_{\mathcal{C}}} & & N
\otimes_{A} {\mathcal{C}} }
\end{equation*}
The set of $\mathcal{C}$-colinear maps from $M$ to $N$ is denoted by $%
\mathrm{Hom}^{\mathcal{C}}(M,N).$ For a right $\mathcal{C}$-comodule $N,$ we
call a right $\mathcal{C}$-comodule $(K,\varrho _{K})$ with $K_{A}\subset
N_{A}$ a $\mathcal{C}$\emph{-subcomodule} if the embedding $K\overset{\iota
_{K}}{\hookrightarrow }N$ is right $\mathcal{C}$-colinear. The category of
counital right $\mathcal{C}$-comodules and right $\mathcal{C}$-colinear maps
is denoted by $\mathbb{M}^{\mathcal{C}}.$ The left $\mathcal{C}$-comodules
and the \emph{left }$C$\emph{-colinear} morphisms are defined analogously.
For left $\mathcal{C}$-comodules $M,N$ we denote with $^{\mathcal{C}}\mathrm{%
Hom}(M,N)$ the set of \emph{left }$\mathcal{C}$\emph{-colinear maps} from $M$
to $N.$ The category of counital\emph{\ left }$\mathcal{C}$\emph{-comodules}
and \emph{left }$\mathcal{C}$\emph{-colinear} morphisms will be denoted with 
$^{\mathcal{C}}\mathbb{M}.$ If $_{A}\mathcal{C}$ (respectively $\mathcal{C}%
_{A})$ is flat, then $\mathbb{M}^{\mathcal{C}}$ (respectively $^{\mathcal{C}}%
\mathbb{M}$) is a Grothendieck category.
\end{punto}

\begin{punto}
\textbf{Measuring }$A$-\textbf{pairings.}\label{mP}\textbf{\ }A \emph{%
measuring left pairing }$P=(\mathcal{T},\mathcal{C}:A)$ consists of an $A$%
-coring $\mathcal{C}$ and an $A$-ring $\mathcal{T}$ with a morphism of $A$%
-rings $\kappa _{P}:\mathcal{T}\rightarrow $ $^{\ast }\mathcal{C}.$ In this
case $\mathcal{C}$ is a right $\mathcal{T}$-module through 
\begin{equation*}
c\leftharpoonup t:=\sum c_{1}<t,c_{2}>\text{ for all }t\in \mathcal{T}\text{
and }c\in \mathcal{C}.
\end{equation*}
Let $(\mathcal{T},\mathcal{C}:A)$ and $(\mathcal{S},\mathcal{D}:B)$ be
measuring left pairings. A (\emph{pure}) \emph{morphism between measuring
pairings} 
\begin{equation*}
(\xi ,\theta :\beta ):(\mathcal{T},\mathcal{C}:A)\rightarrow (\mathcal{S},%
\mathcal{D}:B),
\end{equation*}
consists of an $R$-algebra morphism $\beta :A\rightarrow B,$ a (\emph{pure})
morphism of corings $(\theta :\beta ):(\mathcal{C}:A)\rightarrow (\mathcal{D}%
:B)$ and a morphism of $A$-rings $\xi :\mathcal{S}\rightarrow \mathcal{T}$ $,
$ such that 
\begin{equation*}
<s,\theta (c)>=\beta (<\xi (s),c>)\text{ for all }s\in \mathcal{S}\text{ and 
}c\in \mathcal{C}.
\end{equation*}
If $A=B$ and $\beta $ is not mentioned explicitly, then we mean $\beta
=id_{A}$ and write $(\xi ,\theta ):(\mathcal{T},\mathcal{C})\rightarrow (%
\mathcal{S},\mathcal{D}).$ The category of measuring left pairings and
morphisms described above will be denoted by $\mathcal{P}_{ml}.$
\end{punto}

\begin{punto}
\textbf{The }$\alpha $\textbf{-condition.} We say a measuring left pairing $%
P=(\mathcal{T},\mathcal{C}:A)$ satisfies the $\alpha $\emph{-condition }(or
is a \emph{measuring left }$\alpha $\emph{-pairing}),\emph{\ }if for every
right $A$-module $M$ the following map is injective 
\begin{equation}
\alpha _{M}^{P}:M\otimes _{A}\mathcal{C}\rightarrow \mathrm{Hom}_{-A}(%
\mathcal{T},M),\text{ }\sum m_{i}\otimes _{A}c_{i}\mapsto \lbrack t\mapsto
\sum m_{i}<t,c_{i}>].  \label{alp}
\end{equation}

We say an $A$-coring $\mathcal{C}$ satisfies the \emph{left }$\alpha $\emph{%
-condition}, if the canonical measuring left pairing $(^{\ast }\mathcal{C},%
\mathcal{C}:A)$ satisfies the $\alpha $-condition (equivalently if $_{A}%
\mathcal{C}$ is \emph{locally projective} in the sense of Zimmermann-Huisgen 
\cite[Theorem 2.1]{Z-H76}, \cite[Theorem 3.2]{Gar76}). With $\mathcal{P}%
_{ml}^{\alpha }\subset \mathcal{P}_{ml}$ we denote the \emph{full }%
subcategory of measuring left pairings satisfying the $\alpha $-condition.
The category of \emph{measuring right pairings} $\mathcal{P}_{mr}$ and the
full subcategory of \emph{measuring right }$\alpha $\emph{-pairings} $%
\mathcal{P}_{mr}^{\alpha }\subset \mathcal{P}_{mr}$ are analogously defined.
\end{punto}

\begin{remark}
\label{flat}(\cite[Remark 1.30]{Abu03}) Let $P=(\mathcal{T},\mathcal{C}:A)$
be a measuring left $\alpha $-pairing. Then $_{A}\mathcal{C}$ is flat and $A$%
-cogenerated. If moreover $_{A}\mathcal{C}$ is finitely presented or $_{A}A$
is perfect, then $_{A}\mathcal{C}$ turns out to be projective.
\end{remark}

\begin{punto}
\label{rat-dar}Let\ $P=(\mathcal{T},\mathcal{C}:A)$ be a measuring left $%
\alpha $-pairing and $M$ be a unital right $\mathcal{T}$-module. Consider
the canonical $A$-linear mapping $\rho _{M}:M\rightarrow \mathrm{Hom}_{-A}(%
\mathcal{T},M)$ and put $\mathrm{Rat}^{\mathcal{C}}(M_{\mathcal{T}}):=\rho
_{M}^{-1}(M\otimes _{A}\mathcal{C}).$ We call $M_{\mathcal{T}}$ $\mathcal{C}$%
\emph{-rational}, if $\mathrm{Rat}^{\mathcal{C}}(M_{\mathcal{T}})=M$
(equivalently, if for every $m\in M$ there exists a uniquely determined
element $\sum m_{i}\otimes _{A}c_{i}$ such that $mt=\sum m_{i}<t,c_{i}>$ for
every $t\in \mathcal{T}).$ In this case we get a well-defined $\mathcal{T}$%
-linear map $\varrho _{M}:=(\alpha _{M}^{P})^{-1}\circ \rho
_{M}:M\rightarrow M\otimes _{A}\mathcal{C}.$ The category of $\mathcal{C}$%
-rational right $\mathcal{T}$-modules and right $\mathcal{T}$-linear maps is
denoted by $\mathrm{Rat}^{\mathcal{C}}(\mathbb{M}_{\mathcal{T}}).$ For a
measuring right $\alpha $-pairing $P=(\mathcal{T},\mathcal{C}:A)$ the $%
\mathcal{C}$-rational left $\mathcal{T}$-modules are analogously defined and
we denote by $^{\mathcal{C}}\mathrm{Rat}(_{\mathcal{T}}\mathbb{M})\subseteq $
$_{\mathcal{T}}\mathbb{M}$ the subcategory of $\mathcal{C}$-rational left $%
\mathcal{T}$-modules.
\end{punto}

\begin{proposition}
\emph{(\cite[Proposition 2.8]{Abu03})}\label{comod=sg}\emph{\ }Let $P=(%
\mathcal{T},{\mathcal{C}}:A)$ be a measuring ${}$left pairing. If $_{A}{%
\mathcal{C}}$ is locally projective and $\kappa _{P}(\mathcal{T})\subseteq $ 
$^{\ast }\mathcal{C}$ is dense, then 
\begin{equation}
\begin{tabular}{lllll}
$\mathcal{M}^{{\mathcal{C}}}$ & $\simeq $ & $\mathrm{Rat}^{{\mathcal{C}}}(%
\mathcal{M}_{\mathcal{T}})$ & $=$ & $\sigma \lbrack {\mathcal{C}}_{\mathcal{T%
}}]$ \\ 
& $\simeq $ & $\mathrm{Rat}^{{\mathcal{C}}}(\mathcal{M}_{^{\ast }\mathcal{C}%
})$ & $=$ & $\sigma \lbrack {\mathcal{C}}_{^{\ast }\mathcal{C}}].$%
\end{tabular}
\label{MC=}
\end{equation}
\end{proposition}

\section{Coinduction Functors between Categories of Type $\protect\sigma %
\lbrack M]$}

Our coinduction functor between categories of comodules over corings is
derived from a more general coinduction functor between categories of type $%
\sigma \lbrack M].$ In this section we present this general coinduction
functor and study some of its properties. Throughout this section we fix $R$%
-algebras $A,$ $B$ with a morphism of $R$-algebras $\beta :A\rightarrow B,$
and an $A$-ring $\mathcal{T},$ a $B$-ring $\mathcal{S}$ with a morphism of $%
A $-rings $\xi :\mathcal{S}\rightarrow \mathcal{T}.$

\begin{punto}
\label{IndM}For every right $\mathcal{T}$-module $M,$ the canonical right $B$%
-module $M\otimes _{A}B$ is a right $\mathcal{S}$-module through the right $%
\mathcal{S}$-action 
\begin{equation*}
\sum m_{i}\otimes _{A}b_{i}\leftharpoondown s:=\sum m_{i}\xi (s)\otimes
_{A}b_{i}.
\end{equation*}
It's easy to see that we have in fact a \emph{covariant functor} 
\begin{equation*}
-\otimes _{A}B:\mathbb{M}_{\mathcal{T}}\rightarrow \mathbb{M}_{\mathcal{S}}.
\end{equation*}
\end{punto}

\begin{punto}
For every right $\mathcal{S}$-module $N,$ the canonical right $A$-module $%
\mathrm{Hom}_{-\mathcal{S}}(\mathcal{T}\otimes _{A}B,N)$ is a right $%
\mathcal{T}$-module through 
\begin{equation}
(\varphi \leftharpoonup \widetilde{t})(\sum t_{i}\otimes _{A}b_{i})=\varphi
(\sum \widetilde{t}t_{i}\otimes _{A}b_{i}).  \label{t-act}
\end{equation}
It's easily seen that we have in fact a \emph{covariant }functor 
\begin{equation*}
\mathrm{Hom}_{-\mathcal{S}}(\mathcal{T}\otimes _{A}B,-):\mathbb{M}_{\mathcal{%
S}}\rightarrow \mathbb{M}_{\mathcal{T}}.
\end{equation*}
Moreover $(-\otimes _{A}B,\mathrm{Hom}_{-\mathcal{S}}(\mathcal{T}\otimes
_{A}B,-))$ is an adjoint pair through the functorial canonical isomorphisms 
\begin{equation}
\mathrm{Hom}_{-\mathcal{S}}(M\otimes _{A}B,N)\simeq \mathrm{Hom}_{-\mathcal{S%
}}(M\otimes _{\mathcal{T}}(\mathcal{T}\otimes _{A}B),N)\simeq \mathrm{Hom}_{-%
\mathcal{T}}(M,\mathrm{Hom}_{-\mathcal{S}}(\mathcal{T}\otimes _{A}B,N)).
\label{can-iso}
\end{equation}
\end{punto}

\begin{punto}
\textbf{The coinduction functor }$\mathrm{Coind}_{K}^{L}(-).$\label{IndKL}
Let $K$ be a right $\mathcal{T}$-module and consider the covariant functor 
\begin{equation}
\mathrm{HOM}_{-\mathcal{S}}^{(K)}(\mathcal{T}\otimes _{A}B,-):={\normalsize 
\mathrm{Sp}(\sigma \lbrack }K_{\mathcal{T}}{\normalsize ],-)\circ \mathrm{H}}%
\mathrm{om}_{-\mathcal{S}}{\normalsize (\mathcal{T}}\otimes _{A}B%
{\normalsize ,-):}\text{ }\mathbb{M}_{\mathcal{S}}\rightarrow {\normalsize %
\sigma \lbrack }K_{\mathcal{T}}{\normalsize ].}  \label{Big-Hom-sg}
\end{equation}
For every right $\mathcal{S}$-module $L$ (\ref{Big-Hom-sg}) restricts to the
covariant functor 
\begin{equation}
\mathrm{Coind}_{L}^{K}{\normalsize (-):=\sigma \lbrack }L_{\mathcal{S}}%
{\normalsize ]\rightarrow \sigma \lbrack }K_{\mathcal{T}}{\normalsize ],}%
\text{ }N\mapsto {\normalsize \mathrm{Sp}(\sigma \lbrack }K_{\mathcal{T}}%
{\normalsize ],}\mathrm{Hom}_{-\mathcal{S}}{\normalsize (}\mathcal{T}\otimes
_{A}B{\normalsize ,N))}  \label{spMN}
\end{equation}
defined through the commutativity of the following diagram 
\begin{equation*}
\xymatrix{ {\mathcal M}_{\mathcal S} \ar[rrrrdd]^{{\rm HOM}_{-{\mathcal S}}
^{(K)} ({\mathcal T} \otimes_{A} B,-)} \ar[rrrr]^{ {\rm Hom}_{-{\mathcal S}}
({\mathcal T} \otimes_{A} B,-)} & & & & {\mathcal M}_{\mathcal T} \ar[dd]^{
{\rm Sp} (\sigma[K_{\mathcal T}],-) }\\ & & & & \\ \sigma[L_{\mathcal S}]
\ar@{.>}[rrrr]_{ {\rm Coind} _L ^K (-)} \ar@{^{(}->}[uu] & & & &
\sigma[K_{\mathcal T}] }
\end{equation*}
\qquad 
\end{punto}

\begin{proposition}
\label{main-prop}Let $K$ be a right $\mathcal{T}$-module and $L$ be a right $%
\mathcal{S}$-module. If the induced right $\mathcal{S}$-module $K\otimes
_{A}B$ is $L$-subgenerated, then $-\otimes _{A}B:\mathbb{M}_{\mathcal{T}%
}\rightarrow \mathbb{M}_{\mathcal{S}}$ restricts to a functor 
\begin{equation*}
-\otimes _{A}B:\sigma \lbrack K_{\mathcal{T}}]\rightarrow \sigma \lbrack L_{%
\mathcal{S}}].
\end{equation*}
\end{proposition}

\begin{Beweis}
Assume $K\otimes _{A}B$ to be $L$-subgenerated. Let $M\in \sigma \lbrack K_{%
\mathcal{T}}]$ and $\sum\limits_{i=1}^{q}m_{i}\otimes _{A}b_{i}\in M\otimes
_{A}B$ be an arbitrary element. Since $M_{\mathcal{T}}$ is $K$-subgenerated,
there exists for each $i=1,...,q$ a finite subset $%
\{k_{1}^{(i)},...,k_{n_{i}}^{(i)}\}\subset K$ with $\mathrm{An}_{-\mathcal{T}%
}(\{k_{1}^{(i)},...,k_{n_{i}}^{(i)}\})\subseteq \mathrm{An}_{-\mathcal{T}%
}(m_{i}).$ By assumption $K\otimes _{A}B$ is $L$-subgenerated and so for
every $i=1,...,q$ and $j=1,...,n_{i},$ there exists a finite subset $%
W_{j}^{(i)}\subseteq L,$ such that $\mathrm{An}_{-\mathcal{S}%
}(W_{j}^{(i)})\subseteq \mathrm{An}_{-\mathcal{S}}(\{k_{j}^{(i)}\otimes
_{A}1_{B}\}).$ Consider the finite subset $W:=\bigcup\limits_{i=1}^{q}(%
\bigcup\limits_{j=1}^{n_{i}}W_{j}^{(i)})\subset L.$ If $s\in \mathrm{An}_{-%
\mathcal{S}}(W),$ then obviously $\xi (s)\in \mathrm{An}_{-\mathcal{T}%
}(\{m_{1},...,m_{q}\})$ and we have $\sum_{i=1}^{q}m_{i}\otimes
_{A}b_{i}\leftharpoonup s:=\sum\limits_{i=1}^{q}m_{i}\xi (s)\otimes
_{A}b_{i}=0,$ i.e. $\mathrm{An}_{-\mathcal{S}}(W)\subseteq \mathrm{An}_{-%
\mathcal{S}}(\sum\limits_{i=1}^{q}m_{i}\otimes _{A}b_{i}).$ We conclude then
that $M\otimes _{A}B$ is $L$-subgenerated. We are done, since $\sigma
\lbrack K_{\mathcal{T}}]\subseteq \mathbb{M}_{\mathcal{T}}$ and $\sigma
\lbrack L_{\mathcal{S}}]\subseteq \mathbb{M}_{\mathcal{S}}$ are full
subcategories.$\blacksquare $
\end{Beweis}

\begin{lemma}
Let $K$ be a right $\mathcal{T}$-module and $L$ be a right $\mathcal{S}$%
-module with a right $A$-linear map $\theta :K\rightarrow L.$ If 
\begin{equation}
\xi ((0_{L}:\theta (k)))\subseteq (0_{K}:k)\text{ for every }k\in K,
\label{cond}
\end{equation}
then the induced right $\mathcal{T}$-module $K\otimes _{A}B$ is $L$%
-subgenerated. This is the case in particular, if $\theta $ is injective and 
\begin{equation}
\theta (k\leftharpoonup \xi (s))=\theta (k)\leftharpoonup s\text{ for every }%
s\in \mathcal{S}\text{ and }k\in K.  \label{comp}
\end{equation}
\end{lemma}

\begin{Beweis}
Let $\sum\limits_{i=1}^{n}k_{i}\otimes _{A}b_{i}\in K\otimes _{A}B$ be an
arbitrary element, $V:=\{k_{i}:i=1,...,n\}$ and $W:=\{\theta
(k_{i}):i=1,...,n\}.$ If $s\in \mathrm{An}_{-\mathcal{S}}(W)$ is arbitrary,
then it follows by assumption that $\xi (s)\in \mathrm{An}_{-\mathcal{T}}(V)$
and so 
\begin{equation*}
\sum\limits_{i=1}^{n}k_{i}\otimes _{A}b_{i}\leftharpoonup
s:=\sum\limits_{i=1}^{n}k_{i}\xi (s)\otimes
_{A}b_{i}=\sum\limits_{i=1}^{n}0\otimes _{A}b_{i}=0.
\end{equation*}
So $\mathrm{An}_{-\mathcal{S}}(W)\subseteq (0_{K\otimes
_{A}B}:\sum\limits_{i=1}^{n}k_{i}\otimes _{A}b_{i})$ and consequently $%
\sum\limits_{i=1}^{n}k_{i}\otimes _{A}b_{i}\in \mathrm{Sp}(\sigma \lbrack L_{%
\mathcal{S}}],K\otimes _{A}B).$ Since $\sum\limits_{i=1}^{n}k_{i}\otimes
_{A}b_{i}\in K\otimes _{A}B$ is arbitrary, we conclude that $K\otimes
_{A}B\in \sigma \lbrack L_{\mathcal{S}}].$

Assume now that $\theta $ is injective and that the compatibility condition (%
\ref{comp}) is satisfied. Then we have for arbitrary $k\in K$ and $s\in
(0_{L}:\theta (k)):$%
\begin{equation*}
\theta (k\leftharpoonup \xi (s))=\theta (k)\leftharpoonup s=0,
\end{equation*}
hence $k\leftharpoonup \xi (s)=0$ (since $\theta $ is injective by our
assumption). So condition (\ref{cond}) is satisfied and consequently the
right $\mathcal{S}$-module $K\otimes _{A}B$ is $L$-subgenerated.$%
\blacksquare $
\end{Beweis}

\qquad We are now ready to present the main result in this section

\begin{theorem}
\label{main-sg}Let $K$ be a right $\mathcal{T}$-module, $L$ be a right $%
\mathcal{S}$-module and assume the right $\mathcal{S}$-module $K\otimes _{A}B
$ to be $L$-subgenerated. Then we have an adjoint pair of covariant functors 
$(-\otimes _{A}B,\mathrm{Coind}_{L}^{K}(-)).$
\end{theorem}

\begin{Beweis}
By \cite[45.11]{Wis88} $(\iota (-),\mathrm{Sp}(\sigma \lbrack K_{\mathcal{T}%
}],-))$ is an adjoint pair of covariant functors, and so (\ref{can-iso})
provides us with isomorphisms functorial in $M\in \sigma \lbrack K_{\mathcal{%
T}}]$ and $N\in \sigma \lbrack L_{\mathcal{S}}]:$%
\begin{equation*}
\begin{tabular}{lll}
$\mathrm{Hom}_{-\mathcal{S}}(M\otimes _{A}B,N)$ & $\simeq $ & $\mathrm{Hom}%
_{-\mathcal{S}}((\iota (M)\otimes _{\mathcal{T}}\mathcal{T})\otimes _{A}B),N)
$ \\ 
& $\simeq $ & $\mathrm{Hom}_{-\mathcal{S}}(\iota (M)\otimes _{\mathcal{T}}(%
\mathcal{T}\otimes _{A}B),N)$ \\ 
& $\simeq $ & $\mathrm{Hom}_{-\mathcal{T}}(\iota (M),\mathrm{Hom}_{-\mathcal{%
S}}(\mathcal{T}\otimes _{A}B,N)).$ \\ 
& $\simeq $ & $\mathrm{Hom}_{-\mathcal{T}}(M,\mathrm{Sp}(\sigma \lbrack K_{%
\mathcal{T}}],\mathrm{Hom}_{-\mathcal{S}}(\mathcal{T}\otimes _{A}B,N))).$ \\ 
& $\simeq $ & $\mathrm{Hom}_{-\mathcal{T}}(M,\mathrm{Coind}%
_{L}^{K}(N)).\blacksquare $%
\end{tabular}
\end{equation*}
\end{Beweis}

\section{Coinduction Functors between Categories of Comodules for Measuring $%
\protect\alpha $-Pairings}

In this section we apply our results in the previous section to the special
case of coinduction functors between categories of comodules induced by a
morphism of measuring left $\alpha $-pairing. It turns out that our results
in this setting are direct generalization of the previous results on
coinduction functors between categories of comodules for coalgebras obtained
by the author in his dissertation \cite{Abu01} (see also \cite{Abu}).

\section*{The (Ad-)Induction Functor}

\qquad

\begin{punto}
Let $A$ be an $R$-algebra and $(\mathcal{C},\Delta _{\mathcal{C}%
},\varepsilon _{\mathcal{C}})$ be an $A$-coring, $(M,\varrho _{M})\in 
\mathbb{M}^{\mathcal{C}}$ and $(N,\varrho _{N})\in $ $^{\mathcal{C}}\mathbb{M%
}.$ The \emph{cotensor product} $M\square _{\mathcal{C}}N$ is defined
through the exactness of the following sequence of $R$-modules 
\begin{equation*}
0\longrightarrow M\square _{\mathcal{C}}N\longrightarrow M\otimes _{A}N%
\overset{\varpi }{\longrightarrow }M\otimes _{A}\mathcal{C}\otimes _{A}N,
\end{equation*}
where $\varpi :=\varrho _{M}\otimes _{A}id_{N}-id_{M}\otimes _{A}\varrho
_{N}.$ For more information about the cotensor product of comodules the
reader may consult \cite[Sections 21]{BW03}.
\end{punto}

\begin{punto}
\textbf{The induction functor.}\label{ind-funct} (\cite[Proposition 4.4]
{G-T02})\textbf{\ }Let $(\theta :\beta ):(\mathcal{C}:A)\rightarrow (%
\mathcal{D}:B)$ be a morphism of corings. Then we have a covariant functor,
the so called \emph{induction functor} 
\begin{equation*}
-\otimes _{A}B:\mathbb{M}^{\mathcal{C}}\rightarrow \mathbb{M}^{\mathcal{D}},%
\text{ }M\mapsto M\otimes _{A}B,
\end{equation*}
where the canonical right $B$-module $M\otimes _{A}B$ has a structure of a
right $\mathcal{D}$-comodule through 
\begin{equation}
M\otimes _{A}B\rightarrow (M\otimes _{A}B)\otimes _{B}\mathcal{D}\simeq
M\otimes _{A}\mathcal{D},\text{ }m\otimes _{A}b\mapsto \sum m_{<0>}\otimes
_{A}\theta (m_{<1>})b.  \label{-otB}
\end{equation}
\end{punto}

\begin{punto}
\textbf{The ad-induction functor.}\label{ad-ind}\textbf{\ }(\cite[4.5]{G-T02}%
, \cite[24.7-24.9]{BW03}) Let $(\theta :\beta ):(\mathcal{C}:A)\rightarrow (%
\mathcal{D},B)$ be a morphism of corings. Then the canonical $(B,A)$%
-bimodule $B\otimes _{A}\mathcal{C}$ has a structure of a $(\mathcal{D},%
\mathcal{C})$-bicomodule with the canonical right $\mathcal{C}$-comodule
structure and the left $\mathcal{D}$-comodule structure given by 
\begin{equation*}
\varrho :B\otimes _{A}\mathcal{C}\rightarrow \mathcal{D}\otimes
_{B}(B\otimes _{A}\mathcal{C})\simeq \mathcal{D}\otimes _{A}\mathcal{C},%
\text{ }b\otimes _{A}c\mapsto \sum b\theta (c_{1})\otimes _{A}c_{2}.
\end{equation*}
Assume the morphism of corings $(\theta :\beta ):(\mathcal{C}:A)\rightarrow (%
\mathcal{D},B)$ to be pure. For every right $\mathcal{D}$-comodule $N,$ the
canonical right $A$-submodule $N\square _{\mathcal{D}}(B\otimes _{A}\mathcal{%
C})\subseteq N\otimes _{B}(B\otimes _{A}\mathcal{C})\simeq N\otimes _{A}%
\mathcal{C}$ given by 
\begin{equation*}
N\square _{\mathcal{D}}(B\otimes _{A}\mathcal{C}):=\{\sum n^{i}\otimes
_{A}c^{i}\mid \sum n^{i}\otimes _{B}\theta (c^{i}{}_{1})\otimes
_{A}c^{i}{}_{2}=\sum n^{i}{}_{<0>}\otimes _{B}n^{i}{}_{<1>}\otimes _{A}c^{i}
\end{equation*}
is a right $\mathcal{C}$-comodule by 
\begin{equation*}
\widetilde{\varrho }:N\square _{\mathcal{D}}(B\otimes _{A}\mathcal{C}%
)\rightarrow N\square _{\mathcal{D}}(B\otimes _{A}\mathcal{C})\otimes _{A}%
\mathcal{C}:\sum n^{i}\otimes _{A}c^{i}\mapsto \sum n^{i}\otimes
_{A}c^{i}{}_{1}\otimes _{A}c^{i}{}_{2}.
\end{equation*}
Moreover we have a functor, the so called \emph{ad-induction functor} 
\begin{equation*}
-\square _{\mathcal{D}}(B\otimes _{A}\mathcal{C}):\mathbb{M}^{\mathcal{D}%
}\rightarrow \mathbb{M}^{\mathcal{C}}.
\end{equation*}
\end{punto}

\begin{proposition}
\emph{(\cite[Proposition 4.5]{G-T02})}\label{ad-adj} Let $(\theta :\beta ):(%
\mathcal{C}:A)\rightarrow (\mathcal{D}:B)$ be a pure morphism of corings.
Then $(-\otimes _{A}B,-\square _{\mathcal{D}}(B\otimes _{A}\mathcal{C}))$ is
an adjoint pair of covariant functors.
\end{proposition}

\begin{theorem}
\label{main-meas}Let $P=(\mathcal{T},\mathcal{C}:A)$ and $Q=(\mathcal{S},%
\mathcal{D}:B)$ be measuring left $\alpha $-pairings with a morphism in $%
\mathcal{P}_{ml}^{\alpha }:$%
\begin{equation*}
(\xi ,\theta :\beta ):(\mathcal{T},\mathcal{C}:A)\rightarrow (\mathcal{S},%
\mathcal{D}:B).
\end{equation*}
Then we have a covariant functor 
\begin{equation}
\mathrm{HOM}_{-\mathcal{S}}^{(\mathcal{C})}(\mathcal{T}\otimes _{A}B,-):=%
\mathrm{Rat}^{\mathcal{C}}(-)\circ \mathrm{Hom}_{-\mathcal{S}}{\normalsize (%
\mathcal{T}}\otimes _{A}B{\normalsize ,-):}\text{ }\mathbb{M}_{\mathcal{S}%
}\rightarrow \mathbb{M}^{\mathcal{C}}{\normalsize ,}  \label{HOM-meas}
\end{equation}
which restricts to a functor 
\begin{equation}
\mathrm{Coind}_{\mathcal{D}}^{\mathcal{C}}(-):\mathbb{M}^{\mathcal{D}%
}\rightarrow \mathbb{M}^{\mathcal{C}},\text{ }N\mapsto \mathrm{Rat}^{%
\mathcal{C}}(\mathrm{Hom}_{-\mathcal{S}}(\mathcal{T}\otimes _{A}B,N))
\label{Coind-mesa}
\end{equation}
defined by the commutativity of the following diagram 
\begin{equation*}
\xymatrix{{\Bbb M}_{{\cal D}^{*}} \ar[rrrrdd]^{{\rm HOM}_{-{\mathcal S}}
^{(\cal C)} ({\mathcal T} \otimes_{A} B,-)} \ar[rrrr]^{ {\rm
Hom}_{-{\mathcal S}} ({\mathcal T} \otimes_{A} B,-)} & & & & {\Bbb M}_{{\cal
C}^{*}} \ar[dd]^{ {\rm Rat}^{\cal C} (-)} \\ & & & & \\ {\Bbb M}^{\cal D}
\ar@{.>}[rrrr]_{ {\rm Coind} _D ^C (-)} \ar@{^{(}->}[uu] & & & & {\Bbb
M}^{\cal C} }
\end{equation*}
Moreover $(-\otimes _{A}B,\mathrm{Coind}_{\mathcal{D}}^{\mathcal{C}}(-))$ is
an adjoint pair of covariant functors.
\end{theorem}

\begin{Beweis}
Since $P,Q$ satisfy the $\alpha $-condition, Theorem \ref{comod=sg} implies
that $\mathbb{M}^{\mathcal{C}}\simeq \sigma \lbrack \mathcal{C}_{\mathcal{T}%
}]$ and $\mathbb{M}^{\mathcal{D}}\simeq \sigma \lbrack \mathcal{D}_{\mathcal{%
S}}].$ Moreover the canonical right $B$-module $\mathcal{C}\otimes _{A}B$ is
a right $\mathcal{D}$-comodule by 
\begin{equation*}
\mathcal{C}\otimes _{A}B\mapsto \mathcal{C}\otimes _{A}B\otimes _{B}\mathcal{%
D}\simeq \mathcal{C}\otimes _{A}\mathcal{D},\text{ }c\otimes _{A}b\mapsto
\sum c_{1}\otimes _{A}\theta (c_{2})b,
\end{equation*}
and is hence $\mathcal{D}$-subgenerated as a right $\mathcal{T}$-module. The
result follows now by Theorem \ref{main-sg}.$\blacksquare $
\end{Beweis}

\qquad In the case of a base field, the cotensor functor is equivalent to a
suitable $\mathrm{Hom}$-functor (e.g. \cite[Proposition 3.1]{AW02}).
Analogous to the corresponding result for comodules of coalgebras over
commutative base rings (\cite[Proposition 2.3.12]{Abu01}, \cite[Proposition
2.8]{Abu}) we have

\begin{proposition}
\label{cot=Hom(ot)}Let $B$ be an $R$-algebra, $\mathcal{D}$ a $B$-coring, $%
\mathcal{S}$ a $B$-ring with a morphism of $B$-rings $\kappa :\mathcal{S}%
\rightarrow $ $^{\ast }\mathcal{D}^{\ast }$ and assume the left $B$-pairing $%
Q:=(\mathcal{D},\mathcal{S}:B)\;$to satisfy the $\alpha $-condition.

\begin{enumerate}
\item  Let $(M,\varrho _{M})\in \mathbb{M}^{\mathcal{D}},$ $(N,\varrho
_{N})\in $ $^{\mathcal{D}}\mathbb{M}$ and consider $M\otimes _{B}N$ with the 
\emph{canonical }$\mathcal{S}^{op}$-bimodule structure. Then we have for $%
\sum m_{i}\otimes _{B}n_{i}\in M\otimes _{B}N:$%
\begin{equation*}
\sum m_{i}\otimes _{B}n_{i}\in M\square _{\mathcal{D}}N\Leftrightarrow \sum
s^{op}m_{i}\otimes _{B}n_{i}=\sum m_{i}\otimes _{B}n_{i}s^{op}\text{ for all 
}s^{op}\in \mathcal{S}^{op}.
\end{equation*}

\item  We have an isomorphism of functors for every right $\mathcal{D}$%
-comodule $M$ and every left $\mathcal{D}$-comodule $N:$%
\begin{equation*}
M\square _{\mathcal{D}}N\simeq \text{ }_{\mathcal{S}^{op}}\mathrm{Hom}_{%
\mathcal{S}^{op}}(\mathcal{S}^{op},M\otimes _{B}N).
\end{equation*}
\end{enumerate}
\end{proposition}

\begin{Beweis}
\begin{enumerate}
\item  Let $M\in \mathbb{M}^{\mathcal{D}},$ $N\in $ $^{\mathcal{D}}\mathbb{M}
$ and set $\psi :=\alpha _{M\otimes _{B}N}^{Q}\circ \tau _{(23)}.$ Then we
have: 
\begin{equation*}
\begin{tabular}{llll}
& $\sum m_{i}\otimes _{B}n_{i}\in M\square _{\mathcal{D}}N$ &  &  \\ 
$\Leftrightarrow $ & $\sum m_{i<0>}\otimes _{B}m_{i<1>}\otimes _{B}n_{i}$ & $%
=$ & $\sum m_{i}\otimes _{B}n_{i<-1>}\otimes _{B}n_{i<0>},$ \\ 
$\Leftrightarrow $ & $\psi (\sum m_{i<0>}\otimes _{B}m_{i<1>}\otimes
_{B}n_{i})(s)$ & $=$ & $\psi (\sum m_{i}\otimes _{B}n_{i<-1>}\otimes
_{B}n_{i<0>})(s),$ $\forall $ $s\in \mathcal{S}.$ \\ 
$\Leftrightarrow $ & $\sum m_{i<0>}<s,m_{i<1>}>\otimes _{B}n_{i}$ & $=$ & $%
\sum m_{i}\otimes _{B}<s,n_{i<-1>}>n_{i<0>},$ $\forall $ $s\in \mathcal{S}.$
\\ 
$\Leftrightarrow $ & $\sum s^{op}m_{i}\otimes _{B}n_{i}$ & $=$ & $\sum
m_{i}\otimes _{B}n_{i}s^{op},$ $\forall $ $s^{op}\in \mathcal{S}^{op}.$%
\end{tabular}
\end{equation*}

\item  The isomorphism is given by the morphism 
\begin{equation*}
\gamma _{M,N}:M\square _{\mathcal{D}}N\rightarrow \text{ }_{\mathcal{S}^{op}}%
\mathrm{Hom}_{\mathcal{S}^{op}}(\mathcal{S}^{op},M\otimes _{B}N),\text{ }%
m\otimes _{B}n\mapsto \lbrack s^{op}\mapsto s^{op}m\otimes _{B}n\text{ }%
(=m\otimes _{B}ns^{op})]
\end{equation*}
with inverse $\beta _{M,N}:f\mapsto f(1_{\mathcal{S}^{op}}).$ It is easy to
see that $\gamma _{M,N}$ and $\beta _{M,N}$\ are functorial in $M$ and $%
N.\blacksquare $
\end{enumerate}
\end{Beweis}

As a consequence of Proposition \ref{ad-adj}, Theorem \ref{main-meas} and
Proposition \ref{cot=Hom(ot)} we get

\begin{corollary}
Let $P=(\mathcal{T},\mathcal{C}:A),$ $Q=(\mathcal{S},\mathcal{D}:B)$ be
measuring left $\alpha $-pairings with a morphism in $\mathcal{P}%
_{ml}^{\alpha }:$%
\begin{equation*}
(\xi ,\theta :\beta ):(\mathcal{T},\mathcal{C}:A)\rightarrow (\mathcal{S},%
\mathcal{D}:B).
\end{equation*}
Then we have an isomorphism of functors 
\begin{equation*}
\mathrm{Coind}_{\mathcal{D}}^{\mathcal{C}}(-)\simeq -\square _{\mathcal{D}%
}(B\otimes _{A}\mathcal{C}).
\end{equation*}
If moreover $\kappa _{Q}(\mathcal{S})\subseteq $ $^{\ast }\mathcal{D}^{\ast }
$ then we have an isomorphisms of functors 
\begin{equation*}
\mathrm{Coind}_{\mathcal{D}}^{\mathcal{C}}(-)\simeq -\square _{\mathcal{D}%
}(B\otimes _{A}\mathcal{C})\simeq \mathrm{Hom}_{\mathcal{S}^{op}-\mathcal{S}%
^{op}}(\mathcal{S}^{op},-\otimes _{B}(B\otimes _{A}\mathcal{C})).
\end{equation*}
\end{corollary}

\begin{punto}
\label{th-box}\textbf{Corings over the Same Base Ring.} Let $A$ be an $R$%
-algebra, $\mathcal{C},\mathcal{D}$ be $A$-corings and $\theta :\mathcal{C}%
\rightarrow \mathcal{D}$ be an $A$-coring morphism. Then we have the \emph{%
corestriction functor} 
\begin{equation}
(-)^{\theta }:\mathbb{M}^{\mathcal{C}}\rightarrow \mathbb{M}^{\mathcal{D}},%
\text{ }(M,\varrho _{M})\mapsto (M,(id_{M}\otimes _{A}\theta )\circ \varrho
_{M}).  \label{^th}
\end{equation}
In particular, $\mathcal{C}$ is a $(\mathcal{D},\mathcal{C})$-bicomodule
with the canonical right $\mathcal{C}$-comodule structure and the left $%
\mathcal{D}$-comodule structure by 
\begin{equation*}
^{\mathcal{D}}\varrho _{\mathcal{C}}:\mathcal{C}\overset{\Delta _{\mathcal{C}%
}}{\longrightarrow }\mathcal{C}\otimes _{A}\mathcal{C}\overset{\theta
\otimes _{A}id_{\mathcal{C}}}{\longrightarrow }\mathcal{D}\otimes _{A}%
\mathcal{C}.
\end{equation*}
If $(N,\varrho _{N}^{\mathcal{D}})$ is a right $\mathcal{D}$-comodule and 
\begin{equation*}
\omega _{N,\mathcal{C}}:=\varrho _{N}^{\mathcal{D}}\otimes _{A}id_{\mathcal{C%
}}-id_{N}\otimes _{A}\text{ }^{\mathcal{D}}\varrho _{\mathcal{C}}:N\otimes
_{A}\mathcal{C}\rightarrow N\otimes _{A}\mathcal{D}\otimes _{A}\mathcal{C}
\end{equation*}
is $\mathcal{C}$-pure in $\mathbb{M}_{A}$ (e.g. $_{A}\mathcal{C}$ flat),
then the cotensor product $N\square _{\mathcal{D}}\mathcal{C}$ has a
structure of a right $\mathcal{C}$-comodule by 
\begin{equation*}
N\square _{\mathcal{D}}\mathcal{C}\overset{id_{N}\square _{\mathcal{D}%
}\Delta _{\mathcal{C}}}{\longrightarrow }N\square _{\mathcal{D}}(\mathcal{C}%
\otimes _{A}\mathcal{C})\simeq (N\square _{\mathcal{D}}\mathcal{C})\otimes
_{A}\mathcal{C}
\end{equation*}
and we get the \emph{ad-induction functor} 
\begin{equation*}
-\square _{\mathcal{D}}\mathcal{C}:\mathbb{M}{\normalsize ^{\mathcal{D}%
}\rightarrow \mathbb{M}^{\mathcal{C}},}\text{ }N\mapsto N\square _{\mathcal{D%
}}\mathcal{C}.
\end{equation*}
\end{punto}

\begin{corollary}
Let $A$ be an $R$-algebra, $\mathcal{C},\mathcal{D}$ be $A$-corings and $%
\theta :\mathcal{C}\rightarrow \mathcal{D}$ be a morphism of $A$-corings. If 
$_{A}\mathcal{C}$ is locally projective, then the morphism of $A$-rings $%
^{\ast }\theta :$ $^{\ast }\mathcal{D}\rightarrow $ $^{\ast }\mathcal{C}$
induces a covariant functor 
\begin{equation*}
\mathrm{HOM}_{-^{\ast }\mathcal{D}}^{(\mathcal{C})}(^{\ast }\mathcal{C},-):=%
\mathrm{Rat}^{\mathcal{C}}(-)\circ \mathrm{Hom}_{-^{\ast }\mathcal{D}%
}(^{\ast }\mathcal{C},-):\mathbb{M}_{^{\ast }\mathcal{D}}\rightarrow \mathbb{%
M}^{\mathcal{C}},
\end{equation*}
which restricts to the \emph{coinduction functor} 
\begin{equation*}
\mathrm{Coind}_{\mathcal{D}}^{\mathcal{C}}(-):\mathbb{M}^{\mathcal{D}}%
{\normalsize \rightarrow \mathbb{M}^{\mathcal{C}},}\text{ }N\mapsto \mathrm{%
Rat}^{\mathcal{C}}(\mathrm{Hom}_{-^{\ast }\mathcal{D}}(^{\ast }\mathcal{C}%
,N)).
\end{equation*}
Moreover $((-)^{\theta },\mathrm{Coind}_{\mathcal{D}}^{\mathcal{C}}(-))$ is
an adjoint pair of covariant functors and we have isomorphisms of functors 
\begin{equation*}
\mathrm{Coind}_{\mathcal{D}}^{\mathcal{C}}(-)\simeq -\square _{\mathcal{D}}%
\mathcal{C}\simeq \mathrm{Hom}_{(^{\ast }\mathcal{D}^{\ast })^{op}-(^{\ast }%
\mathcal{D}^{\ast })^{op}}((^{\ast }\mathcal{D}^{\ast })^{op},-\otimes _{A}%
\mathcal{C}).
\end{equation*}
\end{corollary}

\section{The general case}

\qquad In this section we consider the case of coinduction functors induced
by a morphism of corings. Throughout the section, fix $R$-algebras $A,$ $B$
with a morphism of $R$-algebras $\beta :A\rightarrow B,$ an $A$-coring $%
\mathcal{C},$ a $B$-coring $\mathcal{D}$ with a morphism of corings $(\theta
:\beta ):(\mathcal{C}:A)\rightarrow (\mathcal{D}:B)$ and set $^{\#}\mathcal{C%
}:=\mathrm{Hom}_{A-}(\mathcal{C},B).$

\begin{lemma}
\label{Hom(C,B)}

\begin{enumerate}
\item  $^{\#}\mathcal{C}$ is a $(^{\ast }\mathcal{C},^{\ast }\mathcal{D})$%
-bimodule through 
\begin{equation}
(\phi \rightharpoonup h)(c)=\sum h(c_{1}\phi (c_{2}))\text{ and }%
(h\leftharpoonup g)(c)=\sum g(\theta (c_{1})h(c_{2})).  \label{*C-*D}
\end{equation}

\item  If $N$ is a right $^{\ast }\mathcal{D}$-module, then $\mathrm{Hom}%
_{-^{\ast }\mathcal{D}}(^{\#}\mathcal{C},N)$ is a right $^{\ast }\mathcal{C}$%
-module by 
\begin{equation*}
(\varphi \leftharpoonup \psi )(h)=\varphi (\psi \rightharpoonup h)\text{ for
all }\psi \in \text{ }^{\ast }\mathcal{C}\mathbf{,}\text{ }\varphi \in 
\mathrm{Hom}_{-^{\ast }\mathcal{D}}(^{\#}\mathcal{C},N)\text{ and }h\in 
\text{ }^{\#}\mathcal{C}.
\end{equation*}
If $N$ is a left $^{\ast }\mathcal{C}$-module, then $\mathrm{Hom}_{^{\ast }%
\mathcal{C}-}(^{\#}\mathcal{C},N)$ is a left $^{\ast }\mathcal{D}$-module by 
\begin{equation*}
(g\rightharpoonup \varphi )(h)=\varphi (h\leftharpoonup g)\text{\emph{\ }for
all }g\in \text{ }^{\ast }\mathcal{D},\text{ }\varphi \in \mathrm{Hom}%
_{^{\ast }\mathcal{C}-}(^{\#}\mathcal{C},N)\text{ and }h\in \text{ }^{\#}%
\mathcal{C}.
\end{equation*}
\end{enumerate}
\end{lemma}

\begin{Beweis}
\begin{enumerate}
\item  For arbitrary $h\in $ $^{\#}\mathcal{C},$ $\phi ,\psi \in $ $^{\ast }%
\mathcal{C},$ $a\in A$ and $c\in \mathcal{C}$ we have 
\begin{equation*}
(\phi \rightharpoonup h)(ac)=\sum h(ac_{1}\phi (c_{2}))=\sum ah(c_{1}\phi
(c_{2}))=a((\phi \rightharpoonup h)(c))
\end{equation*}
and 
\begin{equation*}
\begin{tabular}{lllll}
$((\phi \ast _{l}\psi )\rightharpoonup h)(c)$ & $=$ & $\sum h(c_{1}(\phi
\ast _{l}\psi )(c_{2}))$ & $=$ & $\sum h(c_{1}\psi (c_{21}\phi (c_{22})))$
\\ 
& $=$ & $\sum h(c_{11}\psi (c_{12}\phi (c_{2})))$ & $=$ & $(\psi
\rightharpoonup h)(\sum c_{1}\phi (c_{2}))$ \\ 
& $=$ & $(\phi \rightharpoonup (\psi \rightharpoonup h))(c),$ &  & 
\end{tabular}
\end{equation*}
i.e. $^{\#}\mathcal{C}$ is a left $^{\ast }\mathcal{C}$-module.

For arbitrary $h\in $ $^{\#}\mathcal{C},$ $f,g\in $ $^{\ast }\mathcal{D},$ $%
a\in A$ and $c\in \mathcal{C}$ we have 
\begin{equation*}
\begin{tabular}{lllll}
$(h\leftharpoonup g)(ac)$ & $=$ & $\sum g(\theta (ac_{1})h(c_{2}))$ & $=$ & $%
\sum g(\beta (a)\theta (c_{1})h(c_{2}))$ \\ 
& $=$ & $\sum \beta (a)g(\theta (c_{1})h(c_{2}))$ & $=$ & $\beta (a)\sum
g(\theta (c_{1})h(c_{2}))$ \\ 
& $=$ & $\beta (a)[(h\leftharpoonup g)(c)]$ & $=$ & $a[(h\leftharpoonup
g)(c)]$%
\end{tabular}
\end{equation*}
and 
\begin{equation*}
\begin{tabular}{lllll}
$(h\leftharpoonup f\ast _{l}g)(c)$ & $=$ & $\sum (f\ast _{l}g)(\theta
(c_{1})h(c_{2}))$ & $=$ & $\sum g(\theta (c_{1})_{1}f(\theta
(c_{1})_{2}h(c_{2})))$ \\ 
& $=$ & $\sum g(\theta (c_{11})f(\theta (c_{12})h(c_{2})))$ & $=$ & $\sum
g(\theta (c_{1})f(\theta (c_{21})h(c_{22})))$ \\ 
& $=$ & $\sum g(\theta (c_{1})(h\leftharpoonup f)(c_{2}))$ & $=$ & $%
((h\leftharpoonup f)\leftharpoonup g)(c),$%
\end{tabular}
\end{equation*}

i.e. $^{\#}\mathcal{C}$ is a right $^{\ast }\mathcal{D}$-module.

To show the compatibility between the left $^{\ast }\mathcal{C}$-action and
the right $^{\ast }\mathcal{D}$-action on $^{\#}\mathcal{C}$ pick arbitrary $%
h\in $ $^{\#}\mathcal{C},$ $\psi \in $ $^{\ast }\mathcal{C},$ $g\in $ $%
^{\ast }\mathcal{D}$ and $c\in \mathcal{C}.$ Then we have 
\begin{equation*}
\begin{tabular}{lllll}
$(\psi \rightharpoonup (h\leftharpoonup g))(c)$ & $=$ & $(h\leftharpoonup
g)(\sum c_{1}\psi (c_{2}))$ & $=$ & $g(\sum \theta (c_{11})h(c_{12}\psi
(c_{2})))$ \\ 
& $=$ & $g(\sum \theta (c_{1})h(c_{21}\psi (c_{22})))$ & $=$ & $g(\sum
\theta (c_{1})(\psi \rightharpoonup h)(c_{2}))$ \\ 
& $=$ & $((\psi \rightharpoonup h)\leftharpoonup g)(c).$ &  & 
\end{tabular}
\end{equation*}

\item  Straightforward.$\blacksquare $
\end{enumerate}
\end{Beweis}

\begin{punto}
Assuming $_{A}\mathcal{C}$ to be locally projective we have the functor 
\begin{equation*}
\mathrm{HOM}_{-^{\ast }\mathcal{D}}(^{\#}\mathcal{C},-):=\mathrm{Rat}^{%
\mathcal{C}}(-)\circ \mathrm{Hom}_{-^{\ast }\mathcal{D}}(^{\#}\mathcal{C},-):%
\mathbb{M}_{^{\ast }\mathcal{D}}\rightarrow \mathbb{M}^{\mathcal{C}},
\end{equation*}
which restricts to a covariant coinduction functor 
\begin{equation*}
\mathcal{G}:\mathbb{M}\mathcal{^{\mathcal{D}}}\rightarrow \mathbb{M}^{%
\mathcal{C}},\text{ }N\mapsto \mathrm{Rat}^{\mathcal{C}}(\mathrm{Hom}%
_{-^{\ast }\mathcal{D}}(^{\#}\mathcal{C},N))
\end{equation*}
defined through the commutativity of the following diagram 
\begin{equation*}
\xymatrix{ {\mathcal M}_{^{*} {\mathcal D}} \ar@{.>}[rrrrrdd]^{ {\rm HOM}_{-
^{*} {\mathcal D} } ^{({\mathcal C})} {(^{\#} {\cal C},-)}} \ar[rrrrr]^{
{\rm Hom}_{-^{*} {\mathcal D}} {(^{\#} {\cal C},-)}} & & & & & {\mathcal
M}_{^{*} {\mathcal C}} \ar[dd]^{ {\rm Rat}^{{\mathcal C}} (-)} \\ & & & & \\
{\mathcal M}^{\mathcal D} \ar@{.>}[rrrrr]_{\mathcal G} \ar@{^{(}->}[uu] & &
& & & {\mathcal M}^{\mathcal C} }
\end{equation*}
\end{punto}

\begin{remark}
Without further assumptions, it is not evident that $(-\otimes _{A}B,%
\mathcal{G})$ is an adjoint pair of covariant functors.
\end{remark}

To get a right adjoint to $-\otimes _{A}B:\mathbb{M}^{\mathcal{C}%
}\rightarrow \mathbb{M}^{\mathcal{D}}$ we modify the definition of $\mathcal{%
G}$. Before we introduce the new version of the coinduction functor, some
technical results are to be proved.

\begin{lemma}
\label{&C}Consider the maps 
\begin{equation*}
\beta \circ -:\text{ }^{\ast }\mathcal{C}\rightarrow \text{ }^{\#}\mathcal{C}%
,\text{ }-\circ \text{ }\theta :\text{ }^{\ast }\mathcal{D}\rightarrow \text{
}^{\#}\mathcal{C}
\end{equation*}
and set $_{\#}\mathcal{C}:=\mathrm{\func{Im}}(\beta \circ -)$ and $_{\#}%
\mathcal{D}:=\mathrm{\func{Im}}(-\circ $ $\theta ).$

\begin{enumerate}
\item  $_{\#}\mathcal{C}$ is an $A$-ring with multiplication 
\begin{equation}
(\beta \circ \phi )\ast (\beta \circ \psi ):=\beta \circ (\phi \ast _{l}\psi
),  \label{mu-C}
\end{equation}
unit $\beta \circ \varepsilon _{\mathcal{C}}$ and 
\begin{equation*}
\beta \circ -:\text{ }^{\ast }\mathcal{C}\rightarrow \text{ }_{\#}\mathcal{C}
\end{equation*}
is a morphisms of $A$-rings. Moreover $_{\#}\mathcal{C}=$ $^{\ast }\mathcal{C%
}\rightharpoonup (\beta \circ \varepsilon _{\mathcal{C}}),$ and hence is a
cyclic $^{\ast }\mathcal{C}$-submodule of $^{\#}\mathcal{C}.$

\item  $_{\#}\mathcal{D}$ is an $A$-ring with multiplication 
\begin{equation}
(f\circ \theta )\ast (g\circ \theta ):=(f\ast _{l}g)\circ \theta ,
\label{mu-D}
\end{equation}
unit $\varepsilon _{\mathcal{D}}\circ \theta $ and 
\begin{equation*}
-\circ \text{ }\theta :\text{ }^{\ast }\mathcal{D}\rightarrow \text{ }_{\#}%
\mathcal{D}
\end{equation*}
is morphisms of $A$-rings. Moreover $_{\#}\mathcal{D}=(\varepsilon _{%
\mathcal{D}}\circ \theta )\leftharpoonup $ $^{\ast }\mathcal{D},$ and hence
is a cyclic $^{\ast }\mathcal{D}$-submodule of $^{\#}\mathcal{C}.$

\item  If $_{\#}\mathcal{D}\subseteq $ $_{\#}\mathcal{C},$ then $_{\#}%
\mathcal{C}\subseteq $ $^{\#}\mathcal{C}$ is a $(^{\ast }\mathcal{C},^{\ast }%
\mathcal{D})$-subbimodule. Analogously, if $_{\#}\mathcal{C}\subseteq $ $%
_{\#}\mathcal{D}$ then $_{\#}\mathcal{D}\subseteq $ $^{\#}\mathcal{C}$ is a $%
(^{\ast }\mathcal{C},^{\ast }\mathcal{D})$-subbimodule.

\item  If $_{\#}\mathcal{D}\subseteq $ $_{\#}\mathcal{C}$ and $N$ is a right 
$^{\ast }\mathcal{D}$-module, then $\mathrm{Hom}_{-^{\ast }\mathcal{D}}(_{\#}%
\mathcal{C},N)$ has a right $^{\ast }\mathcal{C}$-module structure by 
\begin{equation*}
(\varphi \leftharpoonup \psi )(h)=\varphi (\psi \rightharpoonup h)\text{ for
all }\psi \in \text{ }^{\ast }\mathcal{C},\text{ }\varphi \in \mathrm{Hom}%
_{-^{\ast }\mathcal{D}}(_{\#}\mathcal{C},N)\text{ and }h\in \text{ }_{\#}%
\mathcal{C}.
\end{equation*}
If $_{\#}\mathcal{C}\subseteq $ $_{\#}\mathcal{D}$ and $N$ is a left $^{\ast
}\mathcal{C}$-module, then $\mathrm{Hom}_{^{\ast }\mathcal{C}-}(_{\#}%
\mathcal{D},N)$ has a left $^{\ast }\mathcal{D}$-module structure by 
\begin{equation*}
(g\rightharpoonup \varphi )(h)=\varphi (h\leftharpoonup g)\text{\emph{\ }for
all }g\in \text{ }^{\ast }\mathcal{D},\text{ }\varphi \in \mathrm{Hom}%
_{^{\ast }\mathcal{C}-}(_{\#}\mathcal{D},N)\text{ and }h\in \text{ }_{\#}%
\mathcal{D}.
\end{equation*}
\end{enumerate}
\end{lemma}

\begin{Beweis}
\begin{enumerate}
\item  It's obvious, that $(_{\#}\mathcal{C},\ast ,\beta \circ \varepsilon _{%
\mathcal{C}})$ is an $A$-ring and that $\beta \circ -$ is a morphism of $A$%
-rings. For all $\phi \in $ $^{\ast }\mathcal{C}$ and $c\in \mathcal{C}$ we
have 
\begin{equation*}
\begin{tabular}{lllll}
$(\beta \circ \phi )(c)$ & $=$ & $(\beta \circ \phi )(\sum \varepsilon _{%
\mathcal{C}}(c_{1})c_{2}))$ & $=$ & $\sum [(\beta \circ \varepsilon _{%
\mathcal{C}})(c_{1})][(\beta \circ \phi )(c_{2})]$ \\ 
& $=$ & $\sum (\beta \circ \varepsilon _{\mathcal{C}})(c_{1}\phi (c_{2}))$ & 
$=$ & $(\phi \rightharpoonup \beta \circ \varepsilon _{\mathcal{C}})(c),$%
\end{tabular}
\end{equation*}
hence $_{\#}\mathcal{C}=$ $^{\ast }\mathcal{C}\rightharpoonup (\beta \circ
\varepsilon _{\mathcal{C}})\subseteq $ $^{\#}\mathcal{C}$ is a cyclic left $%
^{\ast }\mathcal{C}$-submodule.

\item  It's obvious that $(_{\#}\mathcal{D},\ast ,\varepsilon _{\mathcal{D}%
}\circ \theta )$ is an $A$-ring and $-\circ $ $\theta $ is a morphism of $A$%
-rings. Moreover, for every $g\in $ $^{\ast }\mathcal{D}$ we have 
\begin{equation*}
\begin{tabular}{lllll}
$(g\circ \theta )(c)$ & $=$ & $(g\circ \theta )(\sum c_{1}\varepsilon _{%
\mathcal{C}}(c_{2}))$ & $=$ & $g(\sum \theta (c_{1})(\beta \circ \varepsilon
_{\mathcal{C}})(c_{2}))$ \\ 
& $=$ & $g(\sum \theta (c_{1})(\varepsilon _{\mathcal{D}}\circ \theta
)(c_{2}))$ & $=$ & $((\varepsilon _{\mathcal{D}}\circ \theta )\leftharpoonup
g)(c),$%
\end{tabular}
\end{equation*}
hence $_{\#}\mathcal{D}=(\varepsilon _{\mathcal{D}}\circ \theta
)\leftharpoonup $ $^{\ast }\mathcal{D}\subseteq $ $^{\#}\mathcal{C}$ is a
cyclic right $^{\ast }\mathcal{D}$-submodule.

\item  Assume that $_{\#}\mathcal{D}\subseteq $ $_{\#}\mathcal{C}.$ Then we
have 
\begin{equation*}
\begin{tabular}{lllll}
$_{\#}\mathcal{C}\leftharpoonup $ $^{\ast }\mathcal{D}$ & $=$ & $(^{\ast }%
\mathcal{C}\rightharpoonup \beta \circ \varepsilon _{\mathcal{C}%
})\leftharpoonup $ $^{\ast }\mathcal{D}$ & $=$ & $^{\ast }\mathcal{C}%
\rightharpoonup (\beta \circ \varepsilon _{\mathcal{C}}\leftharpoonup $ $%
^{\ast }\mathcal{D})$ \\ 
& $=$ & $^{\ast }\mathcal{C}\rightharpoonup (\varepsilon _{\mathcal{D}}\circ
\theta \leftharpoonup $ $^{\ast }\mathcal{D})$ & $=$ & $^{\ast }\mathcal{C}%
\rightharpoonup $ $_{\#}\mathcal{D}$ \\ 
& $\subseteq $ & $^{\ast }\mathcal{C}\rightharpoonup $ $_{\#}\mathcal{C}$ & $%
=$ & $_{\#}\mathcal{C},$%
\end{tabular}
\end{equation*}
hence $_{\#}\mathcal{C}\subseteq $ $^{\#}\mathcal{C}$ is closed under the
right $^{\ast }\mathcal{D}$-action (i.e. a right $^{\ast }\mathcal{D}$%
-submodule).

Analogously, if $_{\#}\mathcal{C}\subseteq $ $_{\#}\mathcal{D},$ then 
\begin{equation*}
\begin{tabular}{lllll}
$^{\ast }\mathcal{C}\rightharpoonup $ $_{\#}\mathcal{D}$ & $=$ & $^{\ast }%
\mathcal{C}\rightharpoonup (\varepsilon _{\mathcal{D}}\circ \theta
\leftharpoonup $ $^{\ast }\mathcal{D})$ & $=$ & $(^{\ast }\mathcal{C}%
\rightharpoonup \varepsilon _{\mathcal{D}}\circ \theta )\leftharpoonup $ $%
^{\ast }\mathcal{D}$ \\ 
& $=$ & $(^{\ast }\mathcal{C}\rightharpoonup \beta \circ \varepsilon _{%
\mathcal{C}})\leftharpoonup $ $^{\ast }\mathcal{D}$ & $=$ & $_{\#}\mathcal{C}%
\leftharpoonup $ $^{\ast }\mathcal{D}$ \\ 
& $=$ & $_{\#}\mathcal{D}\leftharpoonup $ $^{\ast }\mathcal{D}$ & $=$ & $%
_{\#}\mathcal{D},$%
\end{tabular}
\end{equation*}
hence $_{\#}\mathcal{D}\subseteq $ $^{\#}\mathcal{C}$ is closed under the
left $^{\ast }\mathcal{C}$-action (i.e. a left $^{\ast }\mathcal{C}$%
-submodule).

\item  Straightforward.$\blacksquare $
\end{enumerate}
\end{Beweis}

\begin{lemma}
\begin{enumerate}
\item  $_{\#}\mathcal{C}$ $\cap $ $_{\#}\mathcal{D}$ is an $A$-ring with
either multiplications \emph{(\ref{mu-C})} or \emph{(\ref{mu-D})} and unit $%
\beta \circ \varepsilon _{\mathcal{C}}=\varepsilon _{\mathcal{D}}\circ
\theta .$ If $_{\#}\mathcal{D}\subseteq $ $_{\#}\mathcal{C}$ \emph{(}%
respectively $_{\#}\mathcal{C}\subseteq $ $_{\#}\mathcal{D}$\emph{)}, then $%
-\circ $ $\theta :$ $^{\ast }\mathcal{D}\rightarrow $ $_{\#}\mathcal{C}$ 
\emph{\ (}respectively $\beta \circ -:$ $^{\ast }\mathcal{C}\rightarrow $ $%
_{\#}\mathcal{D}$\emph{)} is a morphism of $A$-rings.

\item  For all $\phi \in $ $^{\ast }\mathcal{C}$ and $g\in $ $^{\ast }%
\mathcal{D}$ we have 
\begin{equation*}
(\beta \circ \phi )\leftharpoonup g=\phi \rightharpoonup (g\circ \theta ).
\end{equation*}
Hence $_{\#}\mathcal{C}\leftharpoonup $ $^{\ast }\mathcal{D}=$ $^{\ast }%
\mathcal{C}\rightharpoonup $ $_{\#}\mathcal{D}$ is a $(^{\ast }\mathcal{C}%
,^{\ast }\mathcal{D})$-subbimodule of $^{\#}\mathcal{C}.$
\end{enumerate}
\end{lemma}

\begin{Beweis}
\begin{enumerate}
\item  For arbitrary $f\circ \theta =\beta \circ \phi $ and $g\circ \theta
=\beta \circ \psi $ in $_{\#}\mathcal{C}\cap $ $_{\#}\mathcal{D}$ we have 
\begin{equation*}
\begin{tabular}{lllll}
$(f\circ \theta \ast g\circ \theta )(c)$ & $=$ & $(f\ast _{l}g)(\theta (c))$
& $=$ & $\sum g(\theta (c)_{1}f(\theta (c)_{2}))$ \\ 
& $=$ & $\sum g(\theta (c_{1})f(\theta (c_{2})))$ & $=$ & $\sum g(\theta
(c_{1})(\beta \circ \phi )(c_{2}))$ \\ 
& $=$ & $(g\circ \theta )(\sum c_{1}\phi (c_{2}))$ & $=$ & $\beta (\sum \psi
(c_{1}\phi (c_{2}))$ \\ 
& $=$ & $(\beta \circ (\phi \ast _{l}\psi ))(c))$ & $=$ & $((\beta \circ
\phi )\ast (\beta \circ \psi ))(c),$%
\end{tabular}
\end{equation*}
i.e. 
\begin{equation}
(f\circ \theta )\ast (g\circ \theta )=(\beta \circ \phi )\ast (\beta \circ
\psi ).  \label{mutC=multD}
\end{equation}
The last statement follows immediately from (\ref{mutC=multD}).

\item  For all $\phi \in $ $^{\ast }\mathcal{C},$ $g\in $ $^{\ast }\mathcal{D%
}$ and $c\in \mathcal{C}$ we have 
\begin{equation*}
((\beta \circ \phi )\leftharpoonup g)(c)=\sum g(\theta (c_{1})(\beta \circ
\phi )(c_{2}))=\sum g(\theta (c_{1}\phi (c_{2})))=(\phi \rightharpoonup
(g\circ \theta ))(c).\blacksquare 
\end{equation*}
\end{enumerate}
\end{Beweis}

\qquad

\begin{definition}
Let $(\mathcal{C}:A),$ $(\mathcal{D}:B)$ be corings, $(\theta :\beta ):(%
\mathcal{C}:A)\rightarrow (\mathcal{D}:B)$ a morphisms of corings, and
consider the maps 
\begin{equation*}
\beta \circ -:\text{ }^{\ast }\mathcal{C}\rightarrow \text{ }^{\#}\mathcal{C}%
,\text{ }-\circ \text{ }\theta :\text{ }^{\ast }\mathcal{D}\rightarrow \text{
}^{\#}\mathcal{C}.
\end{equation*}
Set $_{\#}\mathcal{C}:=\mathrm{\func{Im}}(\beta \circ -)$ and $_{\#}\mathcal{%
D}:=\mathrm{\func{Im}}(-\circ $ $\theta ).$ We say $(\theta :\beta )$ is a 
\emph{compatible morphism of corings}, provided $_{\#}\mathcal{D}\subseteq
_{\#}\mathcal{C}$ (i.e. for every $f\in $ $^{\ast }\mathcal{C},$ there
exists some $g\in $ $^{\ast }\mathcal{D}$ with $\beta \circ f=g\circ \theta $%
).
\end{definition}

\begin{punto}
\textbf{The coinduction functor.} Assume $_{A}\mathcal{C}$ to be locally
projective and the morphism of corings $(\mathcal{C}:A)\rightarrow \mathcal{%
(D}:B)$ to be compatible. Then we have a covariant functor 
\begin{equation*}
\mathrm{HOM}_{-^{\ast }\mathcal{D}}^{(\mathcal{C})}(_{\#}\mathcal{C},-):=%
\mathrm{Rat}^{\mathcal{C}}(-)\circ \mathrm{Hom}_{-^{\ast }\mathcal{D}}(_{\#}%
\mathcal{C},-):\mathbb{M}_{^{\ast }\mathcal{D}}\rightarrow \mathbb{M}^{%
\mathcal{C}},
\end{equation*}
which restricts to the covariant functor 
\begin{equation*}
\mathrm{Coind}_{\mathcal{D}}^{\mathcal{C}}(-):\mathbb{M}\mathcal{^{\mathcal{D%
}}}\rightarrow \mathbb{M}^{\mathcal{C}},\text{ }N\mapsto \mathrm{Rat}^{%
\mathcal{C}}(\mathrm{Hom}_{-^{\ast }\mathcal{D}}(_{\#}\mathcal{C},N))
\end{equation*}
defined through the commutativity of the following diagram: 
\begin{equation*}
\xymatrix{ {\mathcal M}_{^{*} {\mathcal D}} \ar@{.>}[rrrrrdd]^{ {\rm HOM}_{-
^{*} {\mathcal D} } ^{( {\mathcal C} )} {(_{\#}{\mathcal C},-)}}
\ar[rrrrr]^{ {\rm Hom}_{-^{*} {\mathcal D}} {(_{\#}{\mathcal C},-)}} & & & &
& {\mathcal M}_{^{*} {\mathcal C}} \ar[dd]^{ {\rm Rat}^{\mathcal C} (-)} \\
& & & & & \\ {\mathcal M}^{\mathcal D} \ar@{.>}[rrrrr]_{{\rm Coind}
_{\mathcal D} ^{\mathcal C} (-) } \ar@{^{(}->}[uu] & & & & & {\mathcal
M}^{\mathcal C} }
\end{equation*}
\end{punto}

\begin{lemma}
\label{func-hom}

\begin{enumerate}
\item  If $_{A}\mathcal{C}$ is locally projective then we have functorial
morphisms in $M\in \mathbb{M}^{\mathcal{C}}$ and $N\in \mathbb{M}^{\mathcal{D%
}}:$%
\begin{equation}
\begin{tabular}{cccc}
$\Phi _{M,N}:$ & $\mathrm{Hom}^{\mathcal{D}}(M\otimes _{A}B,N)$ & $%
\rightarrow $ & $\mathrm{Hom}^{\mathcal{C}}(M,\mathrm{Coind}_{\mathcal{D}}^{%
\mathcal{C}}(N))$ \\ 
& $\varkappa $ & $\mapsto $ & $m\mapsto \lbrack h\mapsto \varkappa (\sum
m_{<0>}\otimes _{A}h(m_{<1>}))].$%
\end{tabular}
\label{Phi}
\end{equation}

\item  If $_{A}\mathcal{C},$ $_{B}\mathcal{D}$ are locally projective and $%
(\theta :\beta ):(\mathcal{C}:A)\rightarrow (\mathcal{D}:B)$ is compatible,
then we have functorial morphisms 
\begin{equation}
\begin{tabular}{cccc}
$\Psi _{M,N}:$ & $\mathrm{Hom}^{\mathcal{C}}(M,\mathrm{Coind}_{\mathcal{D}}^{%
\mathcal{C}}(N))$ & $\rightarrow $ & $\mathrm{Hom}^{\mathcal{D}}(M\otimes
_{A}B,N)$ \\ 
& $\zeta $ & $\mapsto $ & $[m\otimes _{A}b\mapsto (\zeta (m)(\beta \circ
\varepsilon _{\mathcal{C}}))b]$%
\end{tabular}
\label{Psi}
\end{equation}
\end{enumerate}
\end{lemma}

\begin{Beweis}
\begin{enumerate}
\item  First of all we prove that $\Phi _{M,N}$ is well-defined for all $%
M\in \mathbb{M}^{\mathcal{C}}$ and $N\in \mathbb{M}^{\mathcal{D}}.$ For all $%
\varkappa \in \mathrm{Hom}^{\mathcal{D}}(M\otimes _{A}B,N),$ $m\in M,$ $\phi
\in $ $^{\ast }\mathcal{C}$ and $h\in $ $_{\#}\mathcal{C}$ we have 
\begin{equation*}
\begin{tabular}{lll}
$\lbrack (\Phi _{M,N}(\varkappa )(m))\leftharpoonup \phi )](h)$ & $=$ & $%
[\Phi _{M,N}(\varkappa )(m)](\phi \rightharpoonup h)$ \\ 
& $=$ & $\varkappa (\sum m_{<0>}\otimes _{A}(\phi \rightharpoonup
h)(m_{<1>}))$ \\ 
& $=$ & $\varkappa (\sum m_{<0>}\otimes _{A}h(m_{<1>1}\phi (m_{<1>2})))$ \\ 
& $=$ & $\varkappa (\sum m_{<0><0>}\otimes _{A}h(m_{<0><1>}\phi (m_{<1>})))$
\\ 
& $=$ & $[\Phi _{M,N}(\varkappa )(\sum m_{<0>}\phi (m_{<1>}))](h)$ \\ 
& $=$ & $[\Phi _{M,N}(\varkappa )(m\leftharpoonup \phi )](h).$%
\end{tabular}
\end{equation*}
So $\Phi _{M,N}(\varkappa )\in \mathrm{Hom}_{-^{\ast }\mathcal{C}}(M,\mathrm{%
Hom}_{-^{\ast }\mathcal{D}}(_{\#}\mathcal{C},N))=\mathrm{Hom}^{\mathcal{C}%
}(M,\mathrm{Coind}_{\mathcal{D}}^{\mathcal{C}}(N)),$ i.e. $\Phi _{M,N}$ is
well-defined. It's easy to see that $\Phi _{M,N}$ is functorial in $M$ and $%
N.$

\item  For all $\zeta \in \mathrm{Hom}^{\mathcal{C}}(M,\mathrm{Coind}_{%
\mathcal{D}}^{\mathcal{C}}(N)),$ $m\in M,$ $b\in B$ and $g\in $ $^{\ast }%
\mathcal{D}$ with $bg\circ \theta =\beta \circ \phi \in $ $_{\#}\mathcal{C}$
we have 
\begin{equation*}
\begin{tabular}{lll}
$\Psi _{M,N}(\zeta )[(m\otimes _{A}b)\leftharpoonup g]$ & $=$ & $\Psi
_{M,N}(\zeta )[\sum (m_{<0>}\otimes _{A}g(\theta (m_{<1>})b)]$ \\ 
& $=$ & $\sum [\zeta (m_{<0>})(\beta \circ \varepsilon _{\mathcal{C}%
})][(bg)(\theta (m_{<1>}))]$ \\ 
& $=$ & $\sum [\zeta (m)_{<0>}(\beta \circ \varepsilon _{\mathcal{C}%
})][(bg\circ \theta )(\zeta (m)_{<1>})]$ \\ 
& $=$ & $\sum [\zeta (m)_{<0>}(\beta \circ \varepsilon _{\mathcal{C}%
})][(\beta \circ \phi )(\zeta (m)_{<1>})]$ \\ 
& $=$ & $\sum [\zeta (m)_{<0>}(\beta \circ \varepsilon _{\mathcal{C}})]\phi
(\zeta (m)_{<1>})$ \\ 
& $=$ & $(\zeta (m)\leftharpoonup \phi )(\beta \circ \varepsilon _{\mathcal{C%
}})$ \\ 
& $=$ & $\zeta (m)(\phi \rightharpoonup (\beta \circ \varepsilon _{\mathcal{C%
}}))$ \\ 
& $=$ & $\zeta (m)(\beta \circ \phi )$ \\ 
& $=$ & $\zeta (m)(bg\circ \theta )$ \\ 
& $=$ & $\zeta (m)((\varepsilon _{\mathcal{D}}\circ \theta )\leftharpoonup
bg)$ \\ 
& $=$ & $[\zeta (m)(\varepsilon _{\mathcal{D}}\circ \theta )b]\leftharpoonup
g$ \\ 
& $=$ & $[\Psi _{M,N}(\zeta )(m\otimes _{A}b)]\leftharpoonup g$%
\end{tabular}
\end{equation*}
So $\Psi _{M,N}(\zeta )\in \mathrm{Hom}_{-^{\ast }\mathcal{D}}(M\otimes
_{A}B,N)=\mathrm{Hom}^{\mathcal{D}}(M\otimes _{A}B,N).$ It's easy to see
that $\Psi _{M,N}$ is functorial in $M$ and $N.\blacksquare $
\end{enumerate}
\end{Beweis}

\qquad We are now in a position that allows us to introduce the main theorem
in this section

\begin{theorem}
\label{main-Th}Let $(\mathcal{C}:A),$ $(\mathcal{D}:B)$ be corings with $_{A}%
\mathcal{C},$ $_{B}\mathcal{D}$ locally projective and $(\theta :\beta ):(%
\mathcal{C}:A)\rightarrow (\mathcal{D}:B)$ be a compatible morphism of
corings. Then $(-\otimes _{A}B,\mathrm{Coind}_{\mathcal{D}}^{\mathcal{C}}(-))
$ is a pair of adjoint functors.

\begin{Beweis}
By Lemma \ref{func-hom} it remains to prove that $\Phi _{M,N}$ and $\Psi
_{M,N}$ are inverse isomorphisms for all $M\in \mathbb{M}^{\mathcal{C}}$ and 
$N\in \mathbb{M}^{\mathcal{D}}.$ For all $\varkappa \in \mathrm{Hom}^{%
\mathcal{D}}(M\otimes _{A}B,N),$ $m\in M$ and $b\in B$ we have 
\begin{equation*}
\begin{tabular}{lll}
$\lbrack (\Psi _{M,N}\circ \Phi _{M,N})(\varkappa )](m\otimes _{A}b)$ & $=$
& $[(\Phi _{M,N}(\varkappa )(m))(\beta \circ \varepsilon _{\mathcal{C}})]b$
\\ 
& $=$ & $[\varkappa (\sum m_{<0>}\otimes _{A}(\beta \circ \varepsilon _{%
\mathcal{C}})(m_{<1>}))]b$ \\ 
& $=$ & $[\varkappa (\sum m_{<0>}\varepsilon _{\mathcal{C}}(m_{<1>})\otimes
_{A}1_{B})]b$ \\ 
& $=$ & $[\varkappa (m\otimes _{A}1_{B})]b$ \\ 
& $=$ & $\varkappa (m\otimes _{A}b),$%
\end{tabular}
\end{equation*}
i.e. $\Psi _{M,N}\circ \Phi _{M,N}=id_{\mathrm{Hom}^{\mathcal{D}}(M\otimes
_{A}B,N)}.$

On the other hand, for all $\zeta \in \mathrm{Hom}^{\mathcal{C}}(M,\mathrm{%
Coind}_{\mathcal{D}}^{\mathcal{C}}(N)),$ $m\in M$ and $h=\beta \circ \phi
\in $ $_{\#}\mathcal{C}$ we have 
\begin{equation*}
\begin{tabular}{lll}
$\lbrack ((\Phi _{M,N}\circ \Psi _{M,N})(\zeta ))(m)](\beta \circ \phi )$ & $%
=$ & $\Psi _{M,N}(\zeta )(\sum m_{<0>}\otimes _{A}(\beta \circ \phi
)(m_{<1>}))$ \\ 
& $=$ & $\sum [\zeta (m_{<0>})(\beta \circ \varepsilon _{\mathcal{C}%
})][(\beta \circ \phi )(m_{<1>})]$ \\ 
& $=$ & $\sum [\zeta (m)_{<0>}(\beta \circ \varepsilon _{\mathcal{C}})]\phi
(\zeta (m)_{<1>})$ \\ 
& $=$ & $(\zeta (m)\leftharpoonup \phi )(\beta \circ \varepsilon _{\mathcal{C%
}})$ \\ 
& $=$ & $\zeta (m)(\phi \rightharpoonup \beta \circ \varepsilon _{\mathcal{C}%
})$ \\ 
& $=$ & $\zeta (m)(\beta \circ \phi ),$%
\end{tabular}
\end{equation*}
i.e. $\Phi _{M,N}\circ \Psi _{M,N}=id_{\mathrm{Hom}^{\mathcal{C}}(M,\mathrm{%
Coind}_{\mathcal{D}}^{\mathcal{C}}(N))}.\blacksquare $
\end{Beweis}
\end{theorem}

\begin{corollary}
Let $(\mathcal{C}:A),$ $(\mathcal{D}:B)$ be corings with $_{A}\mathcal{C},$ $%
_{B}\mathcal{D}$ locally projective and $(\theta :\beta ):(\mathcal{C}%
:A)\rightarrow (\mathcal{D}:B)$ be a compatible morphism of corings. Then we
have an isomorphism of covariant functors 
\begin{equation*}
\mathrm{Coind}_{\mathcal{D}}^{\mathcal{C}}(-)\simeq -\square _{\mathcal{D}%
}(B\otimes _{A}\mathcal{C}).
\end{equation*}
\end{corollary}

\bigskip 

\textbf{Acknowledgment: }The author is grateful for the excellent research
facilities and the financial support provided by KFUPM.


\bigskip 

\begin{thebibliography}{BCMZ01}
\bibitem[Abu]{Abu}  J.Y. Abuhlail, \emph{Hopf pairings and induction
functors over rings}, \textbf{to appear in }Journal of Algebra and its
Applications.

\bibitem[Abu03]{Abu03}  J.Y. Abuhlail, \emph{Rational modules for corings,}
Commun. Algebra \textbf{31(12),} 5793-5840 (2003).

\bibitem[Abu01]{Abu01}  J.Y. Abuhlail, \emph{Dualit\"{a}tstheoreme f\"{u}r
Hopf-Algebren \"{u}ber Ringen}, \textbf{Ph.D. Dissertation,} Heinrich-Heine
Universit\"{a}t, D\"{u}sseldorf-Germany (2001).\linebreak
http://www.ulb.uni-duesseldorf.de/diss/mathnat/2001/abuhlail.html (2001).

\bibitem[AW02]{AW02}  L. Abrams and C. Weibel, \emph{Cotensor products for
modules}, Trans. Am. Math. Soc. \textbf{354(6)}, 2173-2185 (2002).

\bibitem[BCMZ01]{BCMZ01}  T. Brzezi\'{n}ski, S. Caenepeel, G. Militaru and
S. Zhu, \emph{Frobenius and Maschke type theorems for Doi-Hopf modules and
entwined modules revisited: a unified approach,} Granja (ed.) et al., Ring
theory and algebraic geometry, Marcel Dekker. Lect. Notes Pure Appl. Math. 
\textbf{221}, 1-31 (2001).

\bibitem[Brz99]{Brz99}  T. Brzezi\'{n}ski, \emph{On modules associated to
coalgebra Galois extensions}, J. Algebra \textbf{215(1)}, 290-317 (1999).

\bibitem[BW03]{BW03}  T. Brzezi\'{n}ski and R. Wisbauer, \emph{Corings and
Comodules, }Lond. Math. Soc. Lec. Not. Ser. \textbf{309, }Cambridge
University Press (2003).

\bibitem[Gar76]{Gar76}  G. Garfinkel, \emph{Universally torsionless and
trace modules}, J. Amer. Math. Soc. \textbf{215}, 119-144 (1976).

\bibitem[G-T02]{G-T02}  J. G\'{o}mez-Torrecillas, \emph{Separable functors
in corings, }Int. J. Math. Math. Sci. \textbf{30(4)}, 203-225 (2002).

\bibitem[Guz89]{Guz89}  F. Guzman, \emph{Cointegrations, relative cohomology
for comodules, and coseparable corings}, J. Algebra \textbf{126}, 211-224
(1989).

\bibitem[Guz85]{Guz85}  F. Guzman, \emph{Cointegrations and Relative
Cohomology for Comodules}, Ph.D. Thesis, Syracuse University USA (1985).

\bibitem[Wis96]{Wis96}  R. Wisbauer, \emph{Modules and Algebras : Bimodule
Structure and Group Actions on Algebras}, 1st ed., Pitman Monographs and
Surveys in Pure and Applied Mathematics, vol.~81, Addison Wesely Longman
Limited (1996).

\bibitem[Wis88]{Wis88}  R. Wisbauer, \emph{Grundlagen der Modul- und
Ringtheorie : Ein Handbuch f\"{u}r Studium und Forschung}, Verlag Reinhard
Fischer, M\"{u}nchen (1988).

\bibitem[Z-H76]{Z-H76}  B. Zimmermann-Huisgen, \emph{Pure submodules of
direct products of free modules}, Math. Ann. \textbf{224}, 233-245 (1976).
\end{thebibliography}
\end{document}